\documentclass[article]{amsart}
\usepackage{enumitem}
\PassOptionsToPackage{final}{graphicx}
\usepackage{nima}
\usepackage{tikz}
\usetikzlibrary{arrows}
\usetikzlibrary{decorations.markings}
\usetikzlibrary{decorations.pathmorphing}
\usepackage{tikz-cd}
\usepackage{mathtools}

\usepackage{amsaddr}

\newtheorem*{note*}{Note}
\newtheorem{note}[thmctr]{Note}

\usetikzlibrary{calc}
\tikzset{
  vertex/.style={circle,minimum size=0.15cm,inner sep=0,fill=black},
  minivertex/.style={circle,minimum size=0.05cm,inner sep=0,fill=black},
  ensquare/.style={rectangle,minimum size=0.3cm,inner sep=0,draw},
  encircle/.style={circle,minimum size=0.25cm,inner sep=0,draw},
  thickeraser/.style={line width=2.4pt, white}
}

\newcounter{apxthmctr}

\theoremstyle{plain}
\newtheorem*{apxthm*}{Theorem}
\newtheorem{apxthm}[apxthmctr]{Theorem}
\newcommand{\apxthmlabel}[1]{\label{apxthm:#1}}
\newcommand{\apxthmref}[1]{\ref{apxthm:#1}}
\newcommand{\Apxthmref}[1]{Theorem~\apxthmref{#1}}

\newtheorem*{apxlem*}{Lemma}
\newtheorem{apxlem}[apxthmctr]{Lemma}
\newcommand{\apxlemlabel}[1]{\label{apxlem:#1}}
\newcommand{\apxlemref}[1]{\ref{apxlem:#1}}
\newcommand{\Apxlemref}[1]{Lemma~\apxlemref{#1}}

\theoremstyle{definition}
\newtheorem*{apxdefn*}{Definition}
\newtheorem{apxdefn}[apxthmctr]{Definition}
\newcommand{\apxdefnlabel}[1]{\label{apxdefn:#1}}

\newtheorem*{apxrmk*}{Remark}

\newcommand{\interior}[1]{\text{interior}(#1)}

\DeclareMathOperator{\Haus}{Haus}

\DeclareMathOperator{\Int}{Int}
\DeclareMathOperator{\Min}{Min}
\DeclareMathOperator{\Area}{Area}
\DeclareMathOperator{\CAT}{CAT}
\DeclareMathOperator{\Th}{Thick}

\newcommand{\defeq}{\vcentcolon=}

\newcommand\ddeg[1]{\textbf{ddeg}(#1)}
\newcommand\lk[1]{\text{link}(#1)}
\renewcommand{\flat}{{\mathbb{E}^2_\square}}

\title{The Quadric Flat Torus Theorem}

\author{Nima Hoda}
\address{Deptartment of Mathematics, Cornell University,
  Ithaca, NY 14853, USA
  \\ \ \\ Instytut Matematyczny,
  Uniwersytet Wroc\l awski\\
  pl.\ Grun\-wal\-dzki 4,
  50--384 Wroc{\l}aw, Poland}
\email{nima@nimahoda.net}
\thanks{Hoda was partially supported by an Natural Sciences and Engineering Research Council of Canada Postdoctoral Fellowship.}

\author{Zachary Munro}
\address{Burnside Hall, McGill University, 805 Sherbrooke St W, Montreal, Quebec H3A~2K6, Canada}
\address{Faculty of Mathematics, Amado Building, Technion, Haifa, 3200003, Israel}
\address{School of Mathematics, Fry Building, University of Bristol, United Kingdom}
\email{zachy.munro@bristol.ac.uk}
\thanks{Munro was partially supported by the National Science Foundation under Award No. DMS--2503331.}

\keywords{quadric group, quadric complex, hereditary modular graph, flat torus theorem}

\subjclass[2010]{20F65, % Geometric group theory
  20F67} % Hyperbolic groups and nonpositively curved groups

\begin{document}

\begin{abstract}
    We prove a flat torus theorem for quadric complexes.  In particular, we show that if a non-cyclic free abelian group $G$ acts metrically properly on a quadric complex $X$, then $G \cong \Z^2$ and $X$ contains a $G$-invariant isometric copy of the regular square tiling of the plane.
    
    Along the way, we also give a complete proof of the fact that any closed surface subgroup in the fundamental group of a combinatorial 2-complex is represented by a combinatorial map from a cellulation of the surface that is locally injective away from vertices.
\end{abstract}

\maketitle

\tableofcontents

\section{Introduction}

Exploiting the nonpositive curvature of a space with a group action in order to understand algebraic properties of the group has a long history, dating back to the work of Dehn on surface groups \cite{Dehn:Surface_word_problem:1912, Dehn:1987}, to small cancellation theory \cite{Greendlinger:1960, Lyndon:2001} and $\CAT(0)$ geometry \cite{Bridson:1999}.  These ideas culminated in the seminal work of Gromov on hyperbolic groups, $\CAT(0)$ groups and asymptotic geometry \cite{Gromov:1987,Gromov:asymptotic_invariants:1993}.  In more recent years, significant attention has been paid to combinatorially nonpositively curved spaces, including $\CAT(0)$ cube complexes \cite{Avann:1961, Nebesky:1971, Gromov:1987, Roller:1998, Gerasimov:1998, Klavzar:1999, Chepoi:2000}, systolic complexes \cite{Soltan:1983, Chepoi:2000, Haglund:2003, Januszkiewicz:2006}, and Helly graphs \cite{Bandelt_Pesch:Dismantling_abs_retracts:1989, Bandelt_Prisner:Clique_graphs_Helly_graphs:1991, chalopin_chepoi_genevois_hirai_osajda:helly:2020}.
This line of work has lead to remarkable applications in other areas of mathematics \cite{Wise:quasiconvex_hierarchy:2021, Agol:2013} while also establishing broad connections between nonpositive curvature in groups and the field of metric graph theory, where many of these classes of spaces had long been studied under different guises.

A recurring pattern in the study of groups acting on nonpositively curved spaces is a Flat Torus Theorem, which in the classical setting of Riemannian manifolds of nonpositive sectional curvature, says that any $\Z^n$ in the fundamental group stabilizes an isometrically embedded copy of Euclidean $n$-space in the universal cover.  This classical theorem has been generalized to $\CAT(0)$ spaces \cite{Bridson:1999}, and combinatorial analogs have been proven for $\CAT(0)$ cube complexes \cite{Woodhouse:Flat_torus:2017}, systolic complexes \cite{Elsner:2009}, etc.

In this article we prove a flat torus theorem for quadric complexes, which are essentially square complexes in which the interior vertices of any minimal area disk diagram have degree at least 4 (i.e. minimal area disk diagrams are $\CAT(0)$ square complexes)  \cite{Hoda:2017}.  Quadric complexes are also precisely the square completions of hereditary modular graphs, which are characterized by the property that every isometrically embedded cycle has length 4  \cite{Bandelt:1988}.  A rich variety of groups act properly cocompactly on quadric complexes, including finitely presented C(4)-T(4) small cancellation groups \cite{Hoda:2017} and some groups that are not cocompactly cubulated \cite{Huang:2016, Haettel:cocompactly_cubulated_artin:2021}.

% The study of groups acting on a space with combinatorial nonpositive curvature is central to geometry group theory. A classical instance in this theme is the theory of small-cancellation complexes, which generalized foundational work of Dehn on surface groups. In particular, complexes satisfying the $C(6)$ and $C(4)$-$T(4)$ small-cancellation conditions are nonpositively curved. The more recently introduced systolic and quadric complexes generalize these classes of complexes, in the sense that $C(6)$ and $C(4)$-$T(4)$ complexes have natural duals which are systolic and quadric, respectively. Quadric complexes are essentially square complexes whose minimal area diagrams are nonpositively curved, in the sense that they are $\CAT(0)$ square complexes.  

% The definitions of quadric and systolic complexes, introduced by group theorists in the search for fruitful notions of nonpositive curvature, were previously explored by metric graph theorists. Both systolic and quadric complexes are determined by their $1$-skeletons, which are precisely the bridged graphs, in the case of systolic complexes, and hereditarily modular graphs, in the case of quadric complexes. Similarly, the $1$-skeletons of $\CAT(0)$ cube complexes are precisely the median graphs.

Our main theorem is as follows.

\begin{mainthm}[The Quadric Flat Torus Theorem]
\label{thm:A}
    Let $G$ be a non-cyclic free abelian group.  Let $G$ act metrically properly on a quadric complex $X$ (e.g. let $G$ act freely on a locally finite quadric complex $X$).  Then $G \isomor \Z^2$ and there is a $G$-invariant flat in $X$.
\end{mainthm}

A \emph{flat} in a quadric complex is an isometrically embedded copy of the standard square tiling of the Euclidean plane. An action of $G$ on a metric space $X$ is \emph{metrically proper} if $\{g\in G: gB\cap B\neq\emptyset\}$ is finite for any metric ball $B\subset X$.   

Our proof of Theorem~\ref{thm:A} largely follows Elsner's proof of the Systolic Flat Torus Theorem \cite{Elsner:2009}. A point of difference is the use of the theory of $\CAT(0)$ square complexes and their dual curves, unavailable in the systolic setting. As an intermediate step in both proofs, it is shown that a locally minimal subcomplex isomorphic to a tiling of the plane is isometrically embedded (\Thmref{localtoflat}). Both of these theorems are combinatorial descendants of the Flat Torus Theorem for $\CAT(0)$ spaces \cite[Theorem~7.1]{Bridson:1999}.

\begin{mainthm}
\label{thm:B}
    Let $G$ have a normal $\Z^2$ subgroup (e.g. let $G$ be virtually $\Z^2$), and suppose $G$ acts metrically properly on a quadric complex $X$.  Then there is a $G$-invariant thickened flat in $X$.
\end{mainthm}

A \emph{thick flat} in a quadric complex is a flat in $X$ ``blown up'' along its vertices.  For a precise description see \Defnref{thickflat} and \Thmref{thickflat}.

During the proofs of Theorems~\ref{thm:A} and \ref{thm:B}, we apply Theorem~\ref{thm:C} below in the case of $S$ being a torus. 

\begin{mainthm}
\label{thm:C}
    Let $X$ be a combinatorial 2-complex and $S$ a closed, orientable surface of nonzero genus. Let $\phi_*:\pi_1 S\to \pi_1 X$ be an embedding. Then there exists a combinatorial cellular structure on $S$ and a combinatorial map $\phi:S\to X$ inducing $\phi_*$.  Moreover, the map $\phi$ can be chosen to be locally injective away from vertices. 	
\end{mainthm}

\subsection{Structure of the paper}

In \Secref{quadric_cplx}, we give background on quadric complexes, proving some lemmas about disk diagrams and geodesics to be used later. In \Secref{isom_flats}, we prove \Thmref{localtoflat}, the main technical result of the paper, which says that a locally area-minimizing map from a square tiling of a plane to a quadric complex is an isometric embedding. In \Secref{thick_flats} we prove \Lemref{equivalentFlats}, which implies that two flats at finite Hausdorff distance are at Hausdorff distance one. In \Secref{flat_torus_theorem}, we prove \Lemref{normalFlat} and \Thmref{flatTorus} which imply Theorem~\ref{thm:A} and Theorem~\ref{thm:B}. In the appendix, we prove \Apxthmref{surfaceMap} and \Apxthmref{surfaceReduction}, which together imply Theorem~\ref{thm:C}. 

% quadric_cplx
% isom_flats
% thick_flats
% flat_torus_theorem
% surface_maps

\subsection{Acknowledgements}

We are grateful to the reviewer whose detailed comments helped greatly improve this article.

\section{Quadric Complexes}
\seclabel{quadric_cplx}

\begin{defn}
    A map $Y\to X$ between cell complexes is \emph{combinatorial} if open cells are mapped homeomorphically onto open cells. A \emph{square complex} $X$ is a 2-complex such that each 2-cell attaching map $\varphi\colon S^1\to X^1$ is combinatorial for some subdivision of $S^1$ as a 4-cycle. Note that such a subdivision is unique.
\end{defn}

\begin{defn}
    The \emph{combinatorial neighborhood} $N(Z)$ of a subcomplex $Z \subset X$ is the union of closed cells intersecting $Z$. If $Z\subset Y\subset X$ are subcomplexes, then $N_Y(Z)\subset X$ denotes the combinatorial neighborhood of $Z$ in $Y$ considered as a subcomplex of $X$. 
\end{defn}

\begin{defn}[Paths and cycles]
    Let $X$ be a complex. A \emph{path} in $X$ is a map $P\to X$ where $P$ is a complex homeomorphic to a closed interval. The path is \emph{trivial} if the interval is trivial. The \emph{length} $|P|$ of a path is the number of $1$-cells in $P$. An \emph{endpoint} of $P \to X$ is the image of an endpoint of $P$. A path is \emph{closed} if its endpoints coincide. A \emph{cycle} in $X$ is a map $C\to X$ where $C$ is a complex homeomorphic to a circle. A path or cycle is \emph{immersed} if it is locally injective. Two cycles (or two paths) $C_1\to X$ and $C_2\to X$ in $X$ are \emph{equivalent} if there exists an isomorphism $C_1\to C_2$ so that $C_1\to X$ factors as $C_1\to C_2\to X$.
\end{defn}

Note that a path or cycle is determined by the sequence of oriented edges that it traverses, so we will sometimes indicate a path by such a sequence.

\begin{defn}[Graph metric]
    The \emph{graph metric} $d:X^0\times X^0\to \mathbb N$ on the $0$-cells of a complex $X$ is defined by 
    $$d(x,y)=\inf\{|P|:P\to X\text{ is a path with endpoints }x,y\}.$$
\end{defn}

\begin{note}
\label{metric}
    The only metric we make use of in this paper is the graph metric. However, we will often blur the distinction between a complex $X$ and its $0$-skeleton $X^0$. For example, when writing ``$G$ acts isometrically on $X$'' it should be understood that the action of $G$ on $X$ restricts to an action by isometries on $X^0$. 
    Denoting the graph metric on a subcomplex $A$ by $d_A$, we say ``$A\subset X$ is isometrically embedded'' to mean the map $(A^0, d_A)\to (X^0,d)$ induced by inclusion is an isometric embedding. 
\end{note}

\begin{defn}[Disk diagrams]
     A \emph{disk diagram} (or just \emph{diagram}) in $X$ is a map $D\to X$ from a finite, contractible 2-complex $D$ with a fixed embedding in the 2-sphere $D\subset S^2$. This embedding induces a cell structure on $S^2$ with a 2-cell $C_\infty\defeq S^2-D$. The \emph{boundary cycle} $\partial_c D\to D$ is the attaching map of $C_\infty$. The \emph{area} $\Area(D)$ of a diagram $D\to X$ is the number of 2-cells in $D$. 
\end{defn}

\begin{figure}[h]
  \centering
  \begin{tikzpicture}[scale=0.85]
    \begin{scope}
      \node[vertex,label={above:$u_0$}] (u0) at (1, 2) {};
      \node[vertex,label={left:$u_1$}] (u1) at (1, 4/3) {};
      \node[vertex,label={left:$u_2$}] (u2) at (1, 2/3) {};
      \node[vertex,label={below:$u_3$}] (u3) at (1, 0) {};
      \draw[thick] (u0) -- (u1);
      \draw[thick] (u1) -- (u2);
      \draw[thick] (u2) -- (u3);
      \draw[thick,postaction={decorate},decoration={markings,mark=at
      position 1/2 with {\arrow{>},\node[label={left:$e_1$}]
        {};}}] (u3) to[out=180,in=270] (0,1) to[out=90,in=180] (u0);
      \draw[thick,postaction={decorate},decoration={markings,mark=at
      position 1/2 with {\arrow{>},\node[label={right:$e_2$}]
        {};}}] (u3) to[out=0,in=270] (2,1) to[out=90,in=0] (u0);
    \end{scope}

    \draw (3,1) -- (3.75,1) node[above] {$f$};
    \draw[->] (3.75,1) -- (4.5,1);
    
    \begin{scope}[xshift=5cm]
      \node[vertex,label={above:$f(u_0)$}] (fu0) at (1, 2) {};
      \node[vertex,label={left:$f(u_1)$}] (fu1) at (1, 4/3) {};
      \node[vertex,label={left:$f(u_2)$}] (fu2) at (1, 2/3) {};
      \node[vertex,label={below:$f(u_3)$}] (fu3) at (1, 0) {};
      \draw[thick] (fu0) -- (fu1);
      \draw[thick] (fu1) -- (fu2);
      \draw[thick] (fu2) -- (fu3);
      \draw[thick,postaction={decorate},decoration={markings,mark=at
      position 1/2 with {\arrow{>},\node[label={right:$e = f(e_i)$}]
        {};}}] (fu3) to[out=0,in=270] (2,1) to[out=90,in=0] (fu0);
    \end{scope}
  \end{tikzpicture}
  \caption{The fold map $f$.  Let $g \colon s_1 \to s_2$ be an
    isomorphism of combinatorial complexes where $s_1$ and $s_2$ are
    squares.  Let $P_1 \subset \bd s_1$ be a combinatorial path of
    length $3$.  The domain of the fold map is
    $s_1 \sqcup s_2 /{\sim}$ where $x \sim g(x)$ for $x \in P_1$.  The
    fold map is the quotient of $s_1 \sqcup s_2 /{\sim}$ identifying
    $[x]_{\sim}$ and $[g(x)]_{\sim}$ for all $x \in s_1$.}
  \figlabel{foldmap}
\end{figure}

\begin{figure}[h]
  \begin{subfigure}{\textwidth}
    \centering
    \begin{tikzpicture}[scale=0.85]
      \begin{scope}
        \node[vertex] (u0) at (1, 2) {};
        \node[vertex] (u1) at (1, 0) {};
        \node[vertex] (v0) at (0, 1) {};
        \node[vertex] (v1) at (1, 1) {};
        \node[vertex] (v2) at (2, 1) {};

        \draw[thick,postaction={decorate},decoration={markings,mark=at
          position 1/2 with {\arrow{>},\node[label={left:$a$}]
            {};}}] (v0) -- (u0);
        \draw[thick] (v1) -- (u0);
        \draw[thick,postaction={decorate},decoration={markings,mark=at
          position 1/2 with {\arrow{>},\node[label={right:$b$}]
            {};}}] (v2) -- (u0);
        \draw[thick,postaction={decorate},decoration={markings,mark=at
          position 1/2 with {\arrow{>},\node[label={left:$c$}]
            {};}}] (u1) -- (v0);
        \draw[thick] (u1) -- (v1);
        \draw[thick,postaction={decorate},decoration={markings,mark=at
          position 1/2 with {\arrow{>},\node[label={right:$d$}]
            {};}}] (u1) -- (v2);
      \end{scope}

      \draw[-implies, double equal sign distance] (3,1) -- (4,1);
      
      \begin{scope}[xshift=5cm]
        \node[vertex] (u0) at (1, 2) {};
        \node[vertex] (u1) at (1, 0) {};
        \node[vertex] (v0) at (0, 1) {};
        \node[vertex] (v2) at (2, 1) {};

        \draw[thick,postaction={decorate},decoration={markings,mark=at
          position 1/2 with {\arrow{>},\node[label={left:$a$}]
            {};}}] (v0) -- (u0);
        \draw[thick,postaction={decorate},decoration={markings,mark=at
          position 1/2 with {\arrow{>},\node[label={right:$b$}]
            {};}}] (v2) -- (u0);
        \draw[thick,postaction={decorate},decoration={markings,mark=at
          position 1/2 with {\arrow{>},\node[label={left:$c$}]
            {};}}] (u1) -- (v0);
        \draw[thick,postaction={decorate},decoration={markings,mark=at
          position 1/2 with {\arrow{>},\node[label={right:$d$}]
            {};}}] (u1) -- (v2);
      \end{scope}
    \end{tikzpicture}
    \caption{}
    \figlabel{repl2}
  \end{subfigure}
  \\
  \begin{subfigure}{\textwidth}
    \centering
    \begin{tikzpicture}[scale=0.85]
      \begin{scope}
        \node[vertex] (u) at (0, 0) {};
        \node[vertex] (b0) at (30:1) {};
        \node[vertex] (a0) at (90:1) {};
        \node[vertex] (b2) at (150:1) {};
        \node[vertex] (a2) at (210:1) {};
        \node[vertex] (b1) at (270:1) {};
        \node[vertex] (a1) at (330:1) {};

        \foreach \v in {a0,a1,a2}
          \draw[thick] (\v) -- (u);

        \draw[thick,postaction={decorate},decoration={markings,mark=at
          position 1/2 with {\arrow{>},\node[label={above:$c$}]
            {};}}] (a0) -- (b0);
        \draw[thick,postaction={decorate},decoration={markings,mark=at
          position 1/2 with {\arrow{>},\node[label={right:$d$}]
            {};}}] (b0) -- (a1);
        \draw[thick,postaction={decorate},decoration={markings,mark=at
          position 1/2 with {\arrow{>},\node[label={below:$e$}]
            {};}}] (a1) -- (b1);
        \draw[thick,postaction={decorate},decoration={markings,mark=at
          position 1/2 with {\arrow{>},\node[label={below:$f$}]
            {};}}] (b1) -- (a2);
        \draw[thick,postaction={decorate},decoration={markings,mark=at
          position 1/2 with {\arrow{>},\node[label={left:$a$}]
            {};}}] (a2) -- (b2);
        \draw[thick,postaction={decorate},decoration={markings,mark=at
          position 1/2 with {\arrow{>},\node[label={above:$b$}]
            {};}}] (b2) -- (a0);
      \end{scope}

      \draw[-implies,double equal sign distance] (1.75,0) -- (2.75,0);
      
      \begin{scope}[xshift=4.5cm]
        \node[vertex] (b0) at (30:1) {};
        \node[vertex] (a0) at (90:1) {};
        \node[vertex] (b2) at (150:1) {};
        \node[vertex] (a2) at (210:1) {};
        \node[vertex] (b1) at (270:1) {};
        \node[vertex] (a1) at (330:1) {};

        \draw[thick,postaction={decorate},decoration={markings,mark=at
          position 1/2 with {\arrow{>},\node[label={above:$c$}]
            {};}}] (a0) -- (b0);
        \draw[thick,postaction={decorate},decoration={markings,mark=at
          position 1/2 with {\arrow{>},\node[label={right:$d$}]
            {};}}] (b0) -- (a1);
        \draw[thick,postaction={decorate},decoration={markings,mark=at
          position 1/2 with {\arrow{>},\node[label={below:$e$}]
            {};}}] (a1) -- (b1);
        \draw[thick,postaction={decorate},decoration={markings,mark=at
          position 1/2 with {\arrow{>},\node[label={below:$f$}]
            {};}}] (b1) -- (a2);
        \draw[thick,postaction={decorate},decoration={markings,mark=at
          position 1/2 with {\arrow{>},\node[label={left:$a$}]
            {};}}] (a2) -- (b2);
        \draw[thick,postaction={decorate},decoration={markings,mark=at
          position 1/2 with {\arrow{>},\node[label={above:$b$}]
            {};}}] (b2) -- (a0);

        \draw[thick] (b2) -- (a1);
      \end{scope}

      \node at (6.25,0) {,};

      \begin{scope}[xshift=8cm]
        \node[vertex] (b0) at (30:1) {};
        \node[vertex] (a0) at (90:1) {};
        \node[vertex] (b2) at (150:1) {};
        \node[vertex] (a2) at (210:1) {};
        \node[vertex] (b1) at (270:1) {};
        \node[vertex] (a1) at (330:1) {};

        \draw[thick,postaction={decorate},decoration={markings,mark=at
          position 1/2 with {\arrow{>},\node[label={above:$c$}]
            {};}}] (a0) -- (b0);
        \draw[thick,postaction={decorate},decoration={markings,mark=at
          position 1/2 with {\arrow{>},\node[label={right:$d$}]
            {};}}] (b0) -- (a1);
        \draw[thick,postaction={decorate},decoration={markings,mark=at
          position 1/2 with {\arrow{>},\node[label={below:$e$}]
            {};}}] (a1) -- (b1);
        \draw[thick,postaction={decorate},decoration={markings,mark=at
          position 1/2 with {\arrow{>},\node[label={below:$f$}]
            {};}}] (b1) -- (a2);
        \draw[thick,postaction={decorate},decoration={markings,mark=at
          position 1/2 with {\arrow{>},\node[label={left:$a$}]
            {};}}] (a2) -- (b2);
        \draw[thick,postaction={decorate},decoration={markings,mark=at
          position 1/2 with {\arrow{>},\node[label={above:$b$}]
            {};}}] (b2) -- (a0);

        \draw[thick] (a0) -- (b1);
      \end{scope}

      \node at (9.75,0) {,};

      \begin{scope}[xshift=11.5cm]
        \node[vertex] (b0) at (30:1) {};
        \node[vertex] (a0) at (90:1) {};
        \node[vertex] (b2) at (150:1) {};
        \node[vertex] (a2) at (210:1) {};
        \node[vertex] (b1) at (270:1) {};
        \node[vertex] (a1) at (330:1) {};

        \draw[thick,postaction={decorate},decoration={markings,mark=at
          position 1/2 with {\arrow{>},\node[label={above:$c$}]
            {};}}] (a0) -- (b0);
        \draw[thick,postaction={decorate},decoration={markings,mark=at
          position 1/2 with {\arrow{>},\node[label={right:$d$}]
            {};}}] (b0) -- (a1);
        \draw[thick,postaction={decorate},decoration={markings,mark=at
          position 1/2 with {\arrow{>},\node[label={below:$e$}]
            {};}}] (a1) -- (b1);
        \draw[thick,postaction={decorate},decoration={markings,mark=at
          position 1/2 with {\arrow{>},\node[label={below:$f$}]
            {};}}] (b1) -- (a2);
        \draw[thick,postaction={decorate},decoration={markings,mark=at
          position 1/2 with {\arrow{>},\node[label={left:$a$}]
            {};}}] (a2) -- (b2);
        \draw[thick,postaction={decorate},decoration={markings,mark=at
          position 1/2 with {\arrow{>},\node[label={above:$b$}]
            {};}}] (b2) -- (a0);

        \draw[thick] (a2) -- (b0);
      \end{scope}
    \end{tikzpicture}
    \caption{}
    \figlabel{repl3}
  \end{subfigure}
  \caption{Replacement rules for quadric complexes.}
  \figlabel{repl}
\end{figure}

\begin{defn}
  \defnlabel{quadric} A \defterm{locally quadric complex} is a square
  complex $X$ satisfying the following conditions.
  \begin{enumerate}
  \item \itmlabel{imm} The attaching map of every square is an immersion.
  \item \itmlabel{piec3} Any diagram in $X$ of the form of the
    domain of the fold map factors through the fold map.  The
    \defterm{fold map} is described in \Figref{foldmap}.
  \item \itmlabel{rep2} For any diagram in $X$ of the form of the
    left-hand side of \Figref{repl2} with immersed boundary, there is
    a diagram in $X$ of the form on the right with the same
    boundary path.
  \item \itmlabel{rep3} For any diagram in $X$ of the form of the
    left-hand side of \Figref{repl3} with immersed boundary, there is
    a diagram in $X$ of one of the forms on the right with the
    same boundary path.
  \end{enumerate}
  A \defterm{quadric complex} is a simply-connected locally quadric
  complex.
\end{defn}

Condition~\defnitmref{quadric}{piec3} implies that no two squares have
the same attaching map.  Conditions \defnitmref{quadric}{rep2} and
\defnitmref{quadric}{rep3} are nonpositive curvature requirements
having important consequences for diagrams in locally quadric
complexes.

Quadric complexes are similar in nature to systolic complexes.  This
is especially apparent in the presentation given by Wise
\cite{Wise:2003}.  Wise also introduces ``generalized
$(p,q)$-complexes'' which encompass systolic complexes as a subclass
of generalized $(3,6)$-complexes and quadric complexes as a subclass
of generalized $(4,4)$-complexes \cite{Wise:2003}.

\begin{lem}[van Kampen's Lemma]
\lemlabel{vanKampen}
    Let $C\to X$ be a nullhomotopic cycle in a cell complex. Then there exists a diagram $D\to X$ such that the composition $\partial_c D \to D \to X$ is equivalent to $C \to X$. 
\end{lem}

\begin{defn}
    Let $X$ be a 2-complex. The \emph{link} $\lk v$ of a vertex $v$ of $X$ is a graph with a vertex for end of an edge incident to $v$. For each corner of a 2-cell between two ends of edges incident to $v$, there is an edge in $\lk v$ joining the corresponding vertices. 
\end{defn}

For a vertex $v$ of a square complex $X$ we use $\deg(v)$ to denote the number of vertices in $\lk v$ and use $\rho(v)$ to denote the number of edges in $\lk v$.

\begin{lem}[Combinatorial Gauss-Bonnet]\lemlabel{gauss_bonnet}
    Let $X$ be a square complex.  The Euler characteristic of $X$ is equal to the sum of the curvatures of the vertices of $X$ where the \emph{curvature} of a vertex $v$ is $\kappa(v) \defeq 1 - \frac{1}{2}\deg(v) + \frac{1}{4}\rho(v)$.
\end{lem}

\subsection{Diagrams in quadric complexes}

In this section we describe some properties of minimal area disk diagrams in quadric complexes.  These are consequences of the fact that such disk diagrams happen to be $\CAT(0)$ square complexes.  For details see \cite[Section~1.2]{Hoda:2017}.

\begin{lem}
\lemlabel{minareaCAT0}
    Let $D\to X$ be a minimal area diagram in a quadric complex. Then $D$ is a $\CAT(0)$ square complex. 
\end{lem}

We will call a diagram that is also a $\CAT(0)$ square complex a \emph{$\CAT(0)$ square complex diagram}. We now collect some lemmas about $\CAT(0)$ square complex diagrams, which by \Lemref{minareaCAT0} will apply to minimal area diagrams in quadric complexes.

\begin{defn}
    A \emph{midcube} of an edge is its midpoint. A \emph{midcube} of a square is the straight, closed line segment joining the midpoints of a pair of opposing edges. Let $D\to X$ be a diagram in a quadric complex. A \emph{dual curve} of $D$ is a nonempty, connected union of midcubes which intersects each square of $D$ in either an empty set or a midcube. 
\end{defn}

\begin{lem}
\lemlabel{dualcurves}
    Let $D$ be a $\CAT(0)$ square complex diagram. The following properties of dual curves hold. 
    \begin{enumerate}
        \item Any dual curve $\gamma$ of $D$ is an embedded segment (possibly a point) such that $D-\gamma$ has two components.
        \item \itmlabel{dualcurvedistance} The distance between 0-cells $x$ and $y$ in $D$ is equal to the number of dual curves $\#(x,y)$ separating $x$ and $y$.
        \item Two dual curves intersect at most once, and the squares of $D$ are in bijective correspondence with the intersection points of pairs of dual curves.
        \item \itmlabel{nodualtriangles} There does not exist a distinct triple of pairwise intersecting dual curves.
    \end{enumerate}
\end{lem}

\begin{cor}\corlabel{isodiametric}
    Let $D$ be a quadric disk diagram with boundary cycle $\partial_c D \to D$.  Then the diameter of $D$ is at most $\frac{1}{2}|\partial_c D|$.
\end{cor}

\begin{proof}
    Each dual curve of $D$ has endpoints on a pair of edges in $\partial_c D$ by \Lemref{dualcurves}(1), and distinct dual curves of $D$ correspond to disjoint pairs of edges of $\partial_c D$. Thus $D$ has at most $\frac{1}{2}|\partial_c D|$ dual curves. By \Lemref{dualcurves}(2), any two 0-cells of $D$ are at distance at most $\frac{1}{2}|\partial_c D|$ in $D$. 
\end{proof}

\begin{lem}\lemlabel{homotope_diagram}
    Let $\widehat D\to X$ be a diagram in a quadric complex. There exists another diagram $D\to X$ such that $D$ is a $\CAT(0)$ square complex, $\partial_c D\to X$ and $\partial_c \widehat D\to X$ are the same cycle, and $D^0\to X$ factors through $\widehat D^0\to X$ via an injection $D^0\to \widehat D^0$.
\end{lem}
\begin{proof}
    By the Cartan-Hadamard Theorem, if $\widehat D$ is not $\CAT(0)$ then it has an interior vertex $v$ of degree less than $4$.  By the proof of \cite[Lemma~1.6]{Hoda:2017}, there is a disk diagram $D_{-} \to X$ for the same cycle $\partial_c \widehat D \to X$ obtained by performing surgery in the combinatorial neighborhood of $v$.  The surgery involves cutting $D$ open along closed curves and filling with new subdiagrams that contain no internal vertices.  Thus $D_{-}^0 \to X$ factors through $\widehat D^0 \to X$.  The resulting diagram $D_{-}$ has fewer squares than $D$ and so, iterating this process, we obtain our desired $\CAT(0)$ disk diagram $D \to X$.
\end{proof}

\begin{lem}\lemlabel{2x2}
    Let $D$ be a $2\times 2$ grid. Let $\widehat D$ be a $\CAT(0)$ square complex diagram and let $f \colon \partial_c D\to \partial_c \widehat D $ be an isomorphism.  Assume that for all vertices $u,v \in \partial_c D$,
    \begin{equation}\label{colipschitz}
        d_D(u,v) \le d_{\widehat D}\bigl(f(u),f(v)\bigr)
    \end{equation}
    where the distances are taken in $D$ and $\widehat D$ after applying the boundary cycle maps.  Then $f$ extends to an isomorphism $D\to \widehat D$. In particular, $\widehat D$ is also a $2\times 2$ grid.
\end{lem}

\begin{proof}
    By \Lemref{dualcurves}(1), every dual curve in a $\CAT(0)$ square complex diagram has its endpoints lying on the boundary.  Thus $\widehat D$ has exactly four dual curves, each of which corresponds to the midpoints of a pair of edges in $\partial_c \widehat D$.  The isomorphism $f:\partial_c D\to \partial_c \widehat D$ allows us to speak of ``corners'' and ``sides'' of $\partial_c \widehat D$.  The inequality \eqref{colipschitz} and \Lemref{dualcurves}(2) imply that each pair of opposite corners of $\partial_c \widehat D$ must be separated by all four dual curves of $\widehat D$. Consequently, each dual curve of $\widehat D$ has its endpoints on opposite sides of $\partial_c \widehat D$. 

    Consider two dual curves $\alpha$, $\beta$ joining opposite sides of $\partial_c \widehat D$. Any other dual curve $\gamma$ intersects both $\alpha$ and $\beta$ so, by \Lemref{dualcurves}(4), $\alpha$ and $\beta$ do not intersect.  That is, dual curves of $\widehat D$ joining opposite sides of $\partial_c \widehat D$ do not cross.  The patterns of dual curves of $D$ and $\widehat D$ are then isomorphic by an isomorphism extending $f$.  This isomorphism of dual curve patterns induces an isomorphism of square complexes $D \to \widehat D$ extending $f$.
\end{proof}

\begin{cor}\corlabel{2x2}
    Let $A\subset X$ be an isometrically embedded $2\times 2$ grid in a quadric complex $X$, and let $D\to X$ be a diagram with $\partial_c D \to X$ equivalent to the inclusion $\partial A \hookrightarrow X$. Then there exists some $0$-cell of $D$ whose image is adjacent the non-corner vertices of $\partial A$. 
\end{cor}

\subsection{Short cycles in quadric complexes}

\begin{lem}
\lemlabel{sixcycles}
    Let $X$ be a quadric complex and let $f \colon C_6 \to X$ be a combinatorial map from the cycle graph $C_6$ of length $6$.  Then, for some antipodal pair of vertices $v, \bar v$ of $C_6$, there is an edge joining $f(v)$ and $f(\bar v)$ in $X$.
\end{lem}
\begin{proof}
    If $f$ is injective then this follows from the graph-theoretic characterization of quadric complexes \cite[Proposition~1.19]{Hoda:2017}.  Otherwise, since $X^1$ is bipartite, some pair of vertices $u$ and $v$ of $C_6$ at distance $2$ satisfy $f(u) = f(v)$.  Then the antipode $\bar v$ of $v$ in $C_6$ is adjacent to $u$ so that $f(\bar v)$ is adjacent to $f(u) = f(v)$.
\end{proof}

\begin{lem}
\lemlabel{filling}
    Any $4$-cycle in a quadric complex is the boundary of a diagram of area at most $1$. Any $6$-cycle is the boundary of a diagram of area at most $2$. 
\end{lem}

\begin{proof}
    Minimal area diagrams are $\CAT(0)$ square complexes by \Lemref{minareaCAT0}. By \Lemref{dualcurves}, there are two dual curves in a minimal area filling of a $4$-cycle. These can cross at most twice so there is at most a single square in the diagram. Similarly, there are exactly three dual curves in a minimal area filling of a $6$-cycle, but by \Lemref{dualcurves}\pitmref{nodualtriangles} not all three can pairwise intersect. Thus there are at most two squares. 
\end{proof}

\subsection{Quadric geodesics}

\begin{defn}
    For integers $m\leq n$, give $[m,n]$ the combinatorial cell structure where the $0$-cells are the integer points and note that a product of $1$-complexes is a square complex. A \emph{ladder diagram} $D\to X$ is a diagram where $D\isom [0,n] \times [0,1]$ with $n\geq 0$. Consider the path $\bigl(\{0\} \times [0,1] \bigr)\cup\bigl([0,n] \times \{0\}\bigr)\cup\bigl(\{n\} \times [0,1]\bigr) \subset \partial D$. If a path $P\to X$ factors as $P\to D\to X$ where $D\to X$ is a ladder diagram and $P\to D$ is the above path, then $D$ is a \emph{shortcut ladder} along $P$. 
\end{defn}

\begin{lem}
\lemlabel{shortcutladder}
    Let $X$ be a quadric complex. A path $\alpha$ in $X$ is not geodesic if and only if there exists a shortcut ladder along some subpath of $\alpha$. 
\end{lem}
\begin{proof}
    Suppose $\alpha$ is not a geodesic. Let $\gamma \colon Q \to X$ be a geodesic with the same endpoints, and let $D \to X$ be a minimal area disk diagram for the concatenation of $\alpha$ and the reverse of $\gamma$.  
    Since $|\gamma| < |\alpha|$, some dual curve of $D$ must start and end on $\alpha$. An innermost dual curve $\sigma$ then gives a shortcut ladder for a subpath of $\alpha$. Indeed, by \Lemref{dualcurves}\pitmref{nodualtriangles}, distinct dual curves crossing $\sigma$ must be disjoint, so $\sigma$ is contained in a shortcut ladder along a subpath of $\alpha$ whose initial and terminal edges contain the endpoints of $\sigma$. The reverse implication is clear because $P \to D$ in the definition of shortcut ladder is already not geodesic.
\end{proof}

\section{Isometric Flats}
\seclabel{isom_flats}

\begin{defn}
    Let $\flat$ denote the standard square tiling of the Euclidean plane, and let $X$ be a quadric complex. 
    A \emph{flat in $X$} is an isometric embedding $\flat\to X$. 
    A map $\flat\to X$ is \emph{locally minimal} if every restriction to the combinatorial neighborhood of an edge is a minimal diagram. 
    A map $\flat\to X$ is \emph{locally isometric} if every restriction to the combinatorial neighborhood of a vertex is an isometric embedding. 
\end{defn}

Later in this section (\Thmref{localtoflat}) we will prove that locally minimal, locally isometric, and isometric embedding are equivalent conditions for a map $\flat\to X$.

\begin{defn}[Knight move]
    Let $F:\flat\to X$ be map to a quadric complex. 
    Identifying the vertices of $\flat$ with $\Z^2$ in the usual way, a \emph{knight move of $F$} is an edge joining vertices of the form $F(x,y), F(x+1,y+2)$ or $F(x,y),F(x+2,y+1)$ or $F(x,y),F(x+2,y-1)$ or $F(x,y),F(x+1,y-2)$.
\end{defn}

A knight move is a particular failure of local isometry.
As it turns out, the absense of knight moves implies a very strong dichotomy (\Thmref{no knight and not minimal implies cone}) for maps $\flat\to X$. 
First, we need two lemmas.

\begin{lem}
\lemlabel{spiders}
    Let $X$ be a quadric complex and let $F \colon \flat\to X$ be a map with no knight moves.
    Consider the combinatorial neighborhood $N(v)\subset \flat$ of a vertex $v$ of $\flat$.
    Let $u_1,s_1,u_2,s_2,u_3,s_3,u_4,s_4$ be the vertices of  $\partial N(v)$ in cyclic order where $s_1,s_2,s_3,s_4$ are the vertices adjacent to $v$ and $u_1,u_2,u_3,u_4$ are the corners of $N(v)$.
    \begin{enumerate}
    \item If a vertex $w\in X$ is adjacent to $F(s_1)$ and $F(s_3)$, then it is also adjacent to $F(s_2)$ and $F(s_4)$.
    \item If a vertex $w\in X$ is adjacent to $F(u_1)$ and $F(u_3)$, then it is also adjacent to $F(u_2)$ and $F(u_4)$. 
    \end{enumerate}
\end{lem}

\begin{proof}
    (1) Consider the embedded $6$-cycle $C = (F(s_1),F(u_2),F(s_2),F(u_3),F(s_3),w)$.  
    By \Lemref{sixcycles}, there must be a pair of antipodal vertices of $C$ that are joined by an edge.  
    Since $F$ has no knight  moves, there do not exist edges joining $F(s_3), F(u_2)$ or $F(s_1), F(u_3)$.
    Thus there exists an edge joining $w$ and $F(s_2)$. 
    A symmetric argument shows $w$ is adjacent to $F(s_4)$. 
    
    (2) Consider the $6$-cycle $C=(F(u_1),F(s_1),F(u_2),F(s_2),F(u_3),w)$. 
    As above, a combination of \Lemref{sixcycles} and the no knight moves assumption implies there exists ane edge joining $w$ and $F(u_2)$.
    A symmetric argument shows $w$ is adjacent to $F(u_4)$.
\end{proof}

\begin{cor}\corlabel{spider_apocalypse}
    Let $X$ be a quadric complex and let $F \colon \flat\to X$ be a map with no knight moves.  Let $v$ be a vertex of $\flat$ and let $u_1$ and $u_3$ be two antipodal corner vertices of $\partial N(v)$.  Let $U$ be the part of the bipartition of the vertex set of $\flat$ that contains both $u_1$ and $u_3$.

    If some vertex $w \in X$ is adjacent to both $F(u_1)$ and $F(u_3)$ then $w$ is adjacent to every vertex of $F(U)$.
\end{cor}
\begin{proof}
    Identify the vertices of $\flat$ with $\Z^2$ in the usual way with $v$ identified with the identity $(0,0)$ and let \[B_n = \bigl\{(x,y) \in \Z^2 \sth |x| + |y| \le n\bigr\}\] and let \[B^{\infty}_n = \bigl\{(x,y) \in \Z^2 \sth \max\{|x|,|y|\} \le n\bigr\}\] for $n \in \N$.
    That is, the sets $B_n$ and $B^{\infty}_n$ are, respectively, the $\ell_1$-ball and the $\ell_{\infty}$-ball of radius $n$ centered at $v$ in the vertex set of $\flat$.
    
    Applying \Lemref{spiders}(2) we see that $w$ is adjacent to the other two corners $u_2$ and $u_4$ of $\partial N(v)$.  Applying \Lemref{spiders}(1) to $u_1$ and $u_2$ (in the combinatorial neighborhood of their common neighbor) we see that $w$ is adjacent to $v$.  Thus $w$ is adjacent to every vertex in $F(B^{\infty}_1 \cap U)$.  Our desired result then follows from the following two claims.

    \begin{claim*}
        If $w$ is adjacent to every vertex in $F(B^{\infty}_n \cap U)$ then $w$ is adjacent to every vertex in $F(B_{2n} \cap U)$
    \end{claim*}
    
    \begin{claim*}
        If $w$ is adjacent to every vertex in $F(B_{2n} \cap U)$ then $w$ is adjacent to every vertex in $F(B^{\infty}_{2n} \cap U)$
    \end{claim*}

    The first claim follows by repeated applications of \Lemref{spiders}(1) and the second claim follows from repeated applications of \Lemref{spiders}(2).
\end{proof}

\begin{lem}\lemlabel{generalized_spurs}
    Let $D \subset \R^2$ be a finite $\CAT(0)$ square complex embedded in the plane.  
    If $D$ is not equal to a single vertex and $D$ has no interior vertices then $\partial D$ has two disjoint cells $A_1$ and $A_2$, each of which is of one of the following types:
    % vertices $A_1 \subset \partial D^0$ and $A_2 \subset \partial D^0$ such that, for $i=1,2$, either
    \begin{enumerate}
    \item a vertex $v\in \partial D$ with curvature $\frac12$, or
    \item a edge $e\subset \partial D$ whose endpoints have curvature $\frac14$. 
    % \item $A_i =\{u\}$ for a vertex $u$ of curvature $\pi$ or
    % \item $A_i = \{v,w\}$ for a pair of adjacent vertices $u,v$ of curvature $\frac{\pi}{2}$.
    \end{enumerate}
\end{lem}
\begin{proof}
    For the purpose of the proof, we refer to an $A_i$ satisfying Condition~(1) as a \emph{spur} of $D$, we refer to an $A_i$ satisfying Condition~(2) as a \emph{square-spur} of $D$, and we refer to an $A_i$ satisfying either condition a \emph{generalized spur} of $D$.
    If $D$ is a single edge then its two vertices are disjoint spurs. 
    If $D$ is a single square then two of its opposing edges are disjoint square-spurs.  
    Thus we can assume that $D$ is not a single cell.  
    Then, since $D$ has no interior vertices, it has a cut vertex or a cut edge.  
    By induction, each closure $D'$ of a component of the complement of this cut cell $c$ satisfies the conclusion of the lemma.  
    Let $A'_1$ and $A'_2$ be disjoint generalized spurs of $D'$.
    
    In the case where $c$ is a vertex, at least one of $A'_1$ or $A'_2$ is also a generalized spur of $D$.
    In the case where $c$ is an edge, it might be that neither $A'_1$ nor $A'_2$ is a generalized spur of $D$ but this can only happen of $A'_1$ and $A'_2$ are both square-spurs intersecting the edge $c$.  
    But then $D'$ is a single square $s$ in which $A'_1$ and $A'_2$ are opposite edges and $c$ is one of the remaining two edges of $s$.  
    Thus the edge of $s$ opposite $c$ is a square-spur of $D$.  
    Therefore in either the case where $c$ is a vertex or $c$ is an edge, each closure of a component of the complement of $c$ contributes a disjoint generalized spur to $D$.
\end{proof}

\begin{thm}\thmlabel{no knight and not minimal implies cone}
    Let $F:\flat\to X$ be a map to a quadric complex. 
    If $F$ has no knight moves, then exactly one of the following two conditions holds:
    \begin{enumerate}
        \item $F$ is locally minimal. 
        \item There exists a vertex $x\in X$ adjacent to every vertex in one part of the bipartition of $F(\flat)$.
    \end{enumerate}
\end{thm}

\begin{proof}
    If Condition~(2) holds then for any edge $e$ of $\flat$, the restriction $F|_{\partial N(e)}$ of $F$ to the boundary of the combinatorial neighborhood $N(e)$ of $e$ has a disk diagram $D \to X$ with at most five squares.  Indeed $x$ is adjacent to the image of every vertex in one part of the bipartition of $\partial N(e)$.  Thus conditions (1) and (2) are mutually exclusive.
    
    Assume now that Condition~(1) does not hold, i.e. $F$ is not locally minimal.  Let $N(e)\subset \flat$ be a combinatorial neighborhood of an edge such that the restriction $F\colon N(e)\to X$ is not a minimal area diagram and let $f\colon D \to X$ be a minimal area disk diagram for $F|_{\partial N(e)}$.
    
    Because $D$ is a $\CAT(0)$ square complex, the combinatorial neighborhood of each interior vertex of $D$ has at least four squares.  Moreover, the combinatorial neighborhoods of two distinct interior vertices can intersect on at most two squares.  Thus, since $|D| < |N(e)| = 6$, there can be at most one interior vertex of $D$.

    We will show that for any vertices $v,w \in \partial N(e)$ at distance $3$ in $\partial N(e)$, their images $F(v),F(w)$ are not joined by an edge in $X$.  In particular, by \Lemref{generalized_spurs}, this excludes the possibility that $D$ has no internal vertex.  Suppose such $v,w \in \partial N(e)$ exist.  
    Since $X^1$ is bipartite and $F$ has no knight moves, $v$ and $w$ are corners of $N(e)$ at distance $3$ in $\partial N(e)$.
    Let $e'$ be the edge of $X$ joining $F(v),F(w)$. 
    There is a length $5$ path $P$ in $N(e)$ such that $F(P)\cup e'$ is a $6$-cycle whose every diagonal is a knight move of $F$, contradicting \Lemref{sixcycles}.

    The problem is reduced to the case where 
    \begin{enumerate}[label=(\Alph*)]
        \item $D$ has a single interior vertex $x$ and
        \item no edge of $X$ joins $F(v),F(w)$ for vertices $v,w\in \partial N(e)$ at distance $3$ in $\partial N(e)$.
    \end{enumerate}
    An immediate consequence of (B) is that $D$ has no disconnecting vertex. 
    We claim $D$ also has no disconnecting edge, i.e. an interior edge joining two vertices of $\partial D$.
    Indeed, since $X^1$ is bipartite, no edge of $D$ joins vertices at even distance in $\partial N(e)$.
    Condition~(B) implies no edge of $D$ joins vertices $v,w\in \partial N(e)$ at distance $3$.
    That leaves only the possibility that a disconnecting edge joins vertices $v,w\in \partial N(e)$ at distance $5$ in $\partial N(e)$.  
    In this case, $D$ is divided into two disk diagrams bounded by $6$-cycles.  But by \Lemref{sixcycles}, this contradicts (B).
    Thus $D$ has five interior edges joining $x$ to every other vertex of $\partial D$.
    Hence $f(x)$ is adjacent to every vertex in one part of the bipartition of $F(\partial N(e))$. 
    It follows from \Corref{spider_apocalypse} that $f(x)$ is adjacent to every vertex in one part of the bipartition of $F(\flat)$. 
\end{proof}

\begin{lem}
\lemlabel{crawlingSpider}
    Let $F:\flat\to X$ be locally isometric with $X$ quadric and give the vertices of $\flat$ their usual labeling with $\Z^2$. Let $[0^\ast,n^\ast]\times[0^\ast, 1^\ast]$ be a copy of $[0,n]\times [0, 1]$. Let $\alpha: [0^\ast,n^\ast] \times [0^\ast, 1^\ast]\to X$ so that the restrictions $[0,n] \times \{0\} \xrightarrow{F} X$ and $[0^\ast,n^\ast] \times \{0^\ast\} \xrightarrow{\alpha} X$ coincide. If $\alpha(0^\ast,1^\ast)$ is adjacent to $F(-1,-1)$ and $F(-1,1)$, then $\alpha(n^\ast,1^\ast)$ is adjacent to $F(n-1,-1)$ and $F(n-1,1)$.
\end{lem}

\begin{figure}[h]
    \centering
    \includegraphics[width=0.4\textwidth]{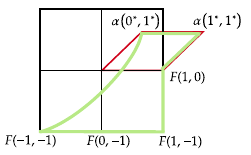}
    \caption{A $6$-cycle to which we apply \Lemref{sixcycles}.}
    \label{fig:spider}
\end{figure}

\begin{proof}
    By induction, it suffices to prove the lemma for the case $n=1$.

    Consider the following $6$-cycle: $F(-1,-1)$, $\alpha(0^\ast,1^\ast)$, $\alpha(1^\ast,1^\ast)$, $F(1,0)$, $F(1,-1)$, $F(0,-1)$. By \Lemref{sixcycles}, there exists an edge joining an antipodal pair of vertices. Since $F$ is locally isometric, an edge cannot join $F(-1,-1)$ and $F(1,0)$, nor $\alpha(0^\ast,1^\ast)$ and $F(1,-1)$. See Figure~\ref{fig:spider}. Thus there must exist an edge joining $\alpha(1^\ast,1^\ast)$ and $F(0,-1)$. A symmetric argument using the $6$-cycle $F(-1,1)$, $\alpha(0^\ast,1^\ast)$, $\alpha(1^\ast,1^\ast)$, $F(1,0)$, $F(1,1)$, $F(0,1)$ shows $\alpha(1^\ast,1^\ast)$ and $F(0,1)$ are joined by an edge.
\end{proof}

\begin{thm}
\thmlabel{localtoflat}
    Let $X$ be a quadric complex. For a map $F:\flat\to X$, the following are equivalent:
    \begin{enumerate}
        \item $F$ is a flat.
        \item $F$ is locally isometric.
        \item $F$ is locally minimal. 
    \end{enumerate}
\end{thm}

\begin{proof}
    $(3)\implies (2)$: Let $D=N(v)$ be the combinatorial neighborhood of some vertex $v$ of $\flat$. Since $D$ is contained in the combinatorial neighborhood of any edge containing $v$, the restriction of $F$ to $D$ is a minimal diagram.  We begin by showing that $F$ is injective on $D^0$. Since all cycles in $X$ are of even length, no two vertices of $D$ at distance one or three can map to the same point of $X$. Since squares in quadric complexes are embedded \cite[Proposition~1.19]{Hoda:2017}, the central vertex $v$ is not mapped to the same vertex as any corner vertex of $D$. If vertices of $\partial D$ at distance two or four map to the same point, then either $\partial D\to X$ factors through a $6$-cycle wedge an edge, or $\partial D\to X$ factors through a wedge of two $4$-cycles. See Figure~\ref{fig:factor1} and Figure~\ref{fig:factor2}. In either case, there is a filling of $\partial D$ with at most two squares by \Lemref{filling}, contradicting minimality. Thus $F$ is injective on $D^0$.

    \begin{figure}[h]
        \centering
        \includegraphics[width=0.5\textwidth]{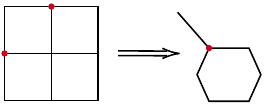}
        \caption{Factoring through a $6$-cycle wedge an edge.}
        \label{fig:factor1}
    \end{figure}

    \begin{figure}[h]
        \centering
        \includegraphics[width=0.5\textwidth]{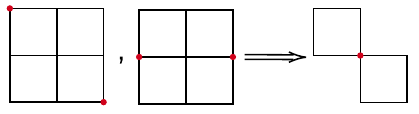}
        \caption{Factoring through a wedge of two $4$-cycles.}
        \label{fig:factor2}
    \end{figure}

    We now show that $F$ is an isometric embedding. Since all cycles in $X$ are of even length, we have only to rule out two cases: vertices at distance three in $D$ map to vertices at distance one in $X$, or vertices at distance four map to vertices at distance two. 

    \begin{figure}[h]
        \centering
        \includegraphics[width=0.5\textwidth]{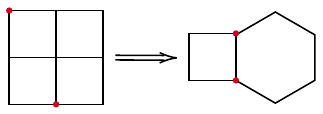}
        \caption{Distance three vertices map to distance one.}
        \label{fig:factor3}
    \end{figure}
    
    In the first case, $\partial D\to X$ factors through a $4$-cycle and $6$-cycle identified along an edge $C_4\cup_e C_6\to X$. See Figure~\ref{fig:factor3}. By \Lemref{filling}, we can fill $\partial D$ with at most three squares, a contradiction. 

    To address the second case, let $u_1$, $u_2$, $u_3$, and $u_4$ be the four corners of $D$, and suppose that there exists a length two path $P\to X$ joining $F(u_1)$ and $F(u_3)$, where $d_D(u_1,u_3)=4$. Let $w$ be the image of the midpoint of $P$. See Figure~\ref{fig:factor4}. There are two $6$-cycles in $F(\partial D)\cup P$, one containing $F(u_1)$, $F(u_2)$, $F(u_3)$, and $w$; and the other containing  $F(u_1)$, $F(u_4)$, $F(u_3)$, and $w$. By \Lemref{sixcycles}, there must exist diagonals in these two $6$-cycles. Since we have already ruled out distance three vertices of $D$ mapping to adjacent vertices in $X$, we must have that $w$ is adjacent to each of $F(u_1)$, $F(u_2)$, $F(u_3)$, and $F(u_4)$.  

    \begin{figure}[h]
        \centering
        \includegraphics[width=0.5\textwidth]{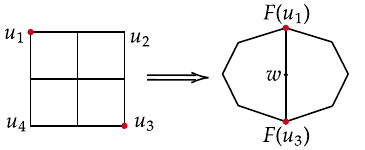}
        \caption{A length two path joins $F(u_1)$ and $F(u_3)$.}
        \label{fig:factor4}
    \end{figure}

    Let $e$ be the edge joining $v$ and the common neighbor $m$ of $u_1$ and $u_2$. The corners $u_3$ and $u_4$ of $D$ are also corners of $N(e)\subset \flat$. The boundary $\partial N(e)$ can be partitioned into subpaths $P_1$, $P_2$, $P_3$, and $P_4$ where $P_i$ joins $u_{i}$ to $u_{i+1}$ for $i\in \Z/4\Z$. Note that $P_1$ has length four, whereas the remaining $P_i$ are length two sides of $D$. The vertex $w$ is adjacent to the endpoints of each $F(P_i)$. Thus there exists a filling of $F(\partial N(e))$ in $X$ by a $6$-cycle and three $4$-cycles. By \Lemref{filling}, $\partial N(e)$ can be filled with at most fives squares, a contradiction. 

    $(2)\implies (1)$: Suppose there exist $p,q\in \flat$ so that $d_X(F(p),F(q))<d_\flat (p,q)$. Since $X^1$ is bipartite, it has no odd cycles and so $d_X(F(p),F(q))$ and $d_\flat (p,q)$ have the same parity. Since $F$ is a local isometry, we have $d_\flat(p,q) \ge 3$. 
    Consider first the case where $d_\flat(p,q)=3$ and  $d_X(F(p),F(q))=1$. Since $F$ is locally isometric, $p$ and $q$ cannot lie in the combinatorial neighborhood of a vertex. Thus $p$ and $q$ lie on a vertical or horizontal line in $\flat$.  Let $p$, $a_1$, $a_2$, $q$ be the vertices of the unique geodesic between $p$ and $q$ in $\flat$. Since $F(p)$ is adjacent to both $F(a_1)$ and $F(q)$, it is also adjacent to the other two neighbors of $F(a_2)$ in $F(N(a_2))$ by \Lemref{spiders}. But then $F(N(a_1))$ is not an isometric embedding, contradicting that $F$ is locally isometric. 

    Suppose now that $d_\flat(p,q)>3$. By \Lemref{shortcutladder}, for any geodesic joining $p$ to $q$ there exists a shortcut ladder $D\to X$ along a subpath $\gamma$ of the geodesic. Moreover, we can choose the path (and so $\gamma$) to be a horizontal path followed by a vertical path. Up to a symmetry of $\flat$, we can assume the endpoints of $\gamma$ are $(-2,0)$ and $(m,n)$ with $m+2\geq n$.  By local isometry, we have $m \ge 1$. By the existence of the shortcut ladder $D\to X$ and an application of \Lemref{spiders}, there exists a map $\alpha:[0^\ast,m^\ast]\times [0^\ast,1^\ast]\to X$, a restriction of $D\to X$, such that $[0,m]\times \{0\} \xrightarrow{F} X$ and $[0^\ast,m^\ast]\times \{0^\ast\} \xrightarrow{\alpha} X$ coincide, and $\alpha(0^\ast, 1^\ast)$ is adjacent to $F(-1,-1)$ and $F(-1,1)$. By \Lemref{crawlingSpider}, $\alpha(m^\ast,1^\ast)$ is adjacent to $F(m-1,-1)$. If $n>1$, a symmetric argument using \Lemref{spiders} and \Lemref{crawlingSpider} on the vertical part of $\gamma$ shows that $\alpha(m^\ast,1^\ast)$ is adjacent to $F(m+1,1)$. Thus $F$ restricted to $N(m,0)$ is not isometric, contradicting local isometry. If $n=1$, then $\alpha((m-1)^\ast,1^\ast)$ is adjacent to $F(m,1)$ since $\alpha$ is a restriction of $D\to X$. An application of \Lemref{crawlingSpider} to the restriction of $\alpha$ on $[0^\ast,(m-1)^\ast]\times[0^\ast,1^\ast]$ shows $\alpha((m-1)^\ast,1^\ast)$ is adjacent to $F(m-2,-1)$. Thus $F$ restricted to $N(m-1,0)$ is not isometric, again contradicting local isometry. Finally, if $n=0$, then by an application of \Lemref{spiders}, $\alpha((m-2)^\ast,1^\ast)$ is adjacent to $F(m-1,1)$, and an application of \Lemref{crawlingSpider} to a restriction of $\alpha$ shows $\alpha((m-2)^\ast,1^\ast)$ is adjacent to $F(m-3,-1)$. Thus $F$ restricted to $N(m-2,0)$ is not isometric, contradicting local isometry.

    $(1)\implies (3)$: %
    % Suppose there exists a combinatorial neighborhood $N(v)$ of a vertex such that $F$ restricted to $N(v)$ is not a minimal diagram. Since $F$ is an isometric embedding, $F$ restricts to an embedding on $\partial N(v)$. In particular, $F(\partial N(v))$ cannot be filled with fewer than three squares. Indeed, there are 8 edges in two squares, and $F(\partial N(v))$ is an 8-cycle, so every edge of a two square filling must be a boundary edge. Thus there exists a disconnecting vertex. If $F(\partial N(v))$ can be filled with three squares, then one of the three squares has three of its sides on $F(\partial N(v))$. Hence $F$ is not isometric.
    Suppose there exists a combinatorial neighborhood $N(e)$ of an edge such that $F$ restricted to $N(e)$ is not a minimal area diagram. $F$ is an embedding on $N(e)$, so $F(\partial N(e))$ cannot be filled by three or fewer squares. If a square has three of its sides on $F(\partial N(e))$, then $F$ is not isometric on $N(e)$. This rules out a filling of $F(\partial N(e))$ with four squares. If $F(\partial N(e))$ can be filled with five squares, then by the above argument, each of the squares must share exactly two edges with $F(\partial N(e))$. For each square, these two edges must be consecutive edges of $F(\partial N(e))$, otherwise $F$ restricted to $\partial N(e)$ is not isometric.  Any such filling must have a unique interior vertex that is adjacent to every vertex in one part of the proper bipartition of the vertices of $F(\partial N(e))$ contradicting $F$ being isometric on $\partial N(e)$.
\end{proof}

\section{Thick Flats}
\seclabel{thick_flats}

% \begin{lem}
%     Let $X$ be a quadric complex, and let $\partial \mathbb{E}_{2\times 2}$ denote the boundary of a square $2\times 2$ grid with the subspace metric. Any minimal area diagram filling an isometric embedding $\partial \mathbb{E}_{2\times 2}\to X$ is isomorphic to $\mathbb{E}_{2\times 2}$.
% \end{lem}

% \begin{proof}
%     Let $D \to X$ be a minimal area diagram filling $\partial \mathbb{E}_{2\times 2}\to X$.  Then $D$ is a $\CAT(0)$ square complex \cite[Section~1.2]{Hoda:2017}.  Since $\partial \mathbb{E}_{2\times 2}\to X$ is an isometric embedding, \Lemref{dualcurves}\pitmref{dualcurvedistance} implies that each antipodal pair of corner vertices of $\partial \mathbb{E}_{2\times 2}$ is separated by all four dual curves of $D$.  It follows that the two dual curves of $D$ starting on a side of $\partial \mathbb{E}_{2\times 2}$ end on the opposing side of $\partial \mathbb{E}_{2\times 2}$.  Then by \Lemref{dualcurves}\pitmref{nodualtriangles}, the two dual curves of $D$ starting on a side of $\partial \mathbb{E}_{2\times 2}$ do not cross.  Thus the pattern of dual curves of $D$ is precisely that of $\mathbb{E}_{2\times 2}$.
% \end{proof}

\begin{defn}
    Let $F:\flat\to X$ be a flat in $X$, and let $v$ be a vertex of $\flat$. Suppose a vertex $w$ of $X$ is adjacent to the images of the neighbors of $v$. Define $F':\flat\to X$ by $F'(v)=w$ and $F'(x)=F(x)$ for all $x\neq v$. Then $F'$ is a flat, and we say $F'$ differs from $F$ by a \emph{vertex move} on $v$. Two flats $F$, $F'$ are \emph{equivalent} if there exists an automorphism $\Phi:\flat\to\flat$ and sequence of flats $F = F_0, F_1, F_2, \ldots$ such that $F_{i+1}$ differs from $F_i$ by a vertex move, for each $i$, and $F_i \to F'\circ\Phi$ pointwise as $i\to \infty$.
\end{defn}

Implicit in the above definition is the assertion that $F'$ is a flat. Indeed, since $F'(x)=F(x)$ for all $x\neq v$, it suffices to check $F'$ preserves distances between $v$ and other vertices. 
Suppose $d_X(F'(v),F'(x))<d_\flat(v,x)$ for some $x\in \flat$. 
Let $y\in \flat$ be a vertex such that some geodesic from $y$ to $x$ passes through $v$. 
Then $d_X(F(y),F(x)) = d_X(F'(y),F'(x)) \le d_X(F'(y),F'(v))+d_X(F'(v),F'(x)) <d_\flat(y,v)+d_\flat(v,x)=d_\flat(y,x)$, contradicting the fact that $F$ is an isometric embedding.

\begin{lem}
\lemlabel{vertexMovesCombine}
    Let $F$ be a flat.  Let $F'$ be obtained from $F$ by a vertex move on $v' \in \flat$ and let $F''$ be obtained from $F$ by a vertex move on $v'' \in \flat$.  If $v'$ and $v''$ are adjacent then $F'(v')$ and $F''(v'')$ are adjacent. Consequently, a collection of possible vertex moves can be performed in any sequence.
\end{lem}

As a consequence of \Lemref{vertexMovesCombine}, any sequence of vertex moves at distinct vertices can be permuted to obtain the same flat.

\begin{proof}
    If $v',v''\in \flat$ are adjacent, then there exist vertices $x_1,x_2\in N(v')-v'$ and $y_1,y_2\in N(v'')- v''$ such that $x_1y_1v''y_2x_2v'$ forms an embedded $6$-cycle in $\flat$. 
    By \Lemref{sixcycles} the $6$-cycle $F(x_1)F(y_1)F''(v'')F(y_2)F(x_2)F'(v')$ contains a diagonal. 
    Since $F$ is an isometric embedding, $\{F(x_1), F(y_2)\}$ and $\{F(x_2),F(y_1)\}$ cannot be edges of $X$. 
    Thus $F'(v')$ and $F''(v'')$ must be adjacent. 
\end{proof}

% \begin{lem}
% \lemlabel{}
%     Let $D\to X$ be a diagram in a quadric complex (not necessarily homeomorphic to a disc) and let $\bd_pD \to D$ be the boundary path of $D$. There exists a planar complex $\mathcal D\supset D$, homeomorphic to a disc, and a map $\mathcal D^1\to X$ extending $D^1\to X$ such that the restriction of $\mathcal D^1\to X$ to $\bd \mathcal D$ is the same cycle in $X$ as the composition $\bd_pD \to D \to X$.
% \end{lem}

% \begin{proof}
%     If $D$ is not homeomorphic to a disc, then there exists a disconnecting $0$-cell $v$ of $D$. The edges incident to $v$ are cyclically ordered from the embedding $D\subset S^2$. For each pair $e_1$, $e_2$ of consecutive edges incident to $v$ belonging to distinct components of $D-v$, we attach a square by gluing in a corner to $e_1\cup e_2$. The map $D^1\to X$ can be extended to each square by sending the square to the length two path $e_1\cup e_2\to X$. In this way we obtain a planar complex $\mathcal D_+ \supset D$ with a map $\mathcal D_+^1\to X$ extending $D^1\to X$ such that $\mathcal D_+$ has one fewer disconnecting $0$-cell than $D$. Repeating this operation for each disconnecting $0$-cell, we get our desired $\mathcal D\supset D$ and map $\mathcal D^1\to X$.
% \end{proof}

\begin{lem}
\lemlabel{equivalentFlats}
    Suppose that the images $A$, $A'$ of two flats $F$, $F'$ in a quadric complex are at finite Hausdorff distance $h\defeq d_{\Haus}(A,A')<\infty$. Then $F$ and $F'$ are equivalent. In particular, $d_{\Haus}(A, A')\leq 1$. 
\end{lem}

\begin{proof}
    For convenience, we identify the domains of $F$ and $F'$ with $A$ and $A'$, respectively, and view $F$ and $F'$ as inclusion maps.

    Without loss of generality, we have $h\geq 1$. Let $v$ be a $0$-cell of $A$. Let $R\subset \flat$ be a square of side-length $10h$ centered at $v$. We subdivide each side of $R$ into ten segments of length $h$. For each endpoint $x$ in this subdivision choose a geodesic $\xi_x$ of length at most $h$ joining $x$ to $A'$. For each segment $\alpha$ of the subdivision with endpoints $x$ and $y$, let $\xi_\alpha$ be a geodesic of $A'$ joining the endpoints of $\xi_x$ and $\xi_y$ contained in $A'$. Note that $|\xi_\alpha|\leq 3h$ by the triangle inequality, so the cycle $\alpha\cup \xi_x\cup\xi_y\cup\xi_\alpha$ has length at most $6h$. Let $D_\alpha\to X$ be a minimal area diagram filling this cycle.

    By construction, the concatenation $\xi$ of the $\xi_\alpha$ is a closed path in $A'$. Thus there exists a diagram $D\to A'$ with $\partial_c D\to A'$ equivalent to $\xi$. Let $R'=R-\Int(N_{A}(v))$ be the annulus obtained from $R$ by deleting the interior of the combinatorial neighborhood of $v$. Identifying diagrams along the boundary segments $\alpha$ of $R$, the $\xi_x$, and the $\xi_\alpha$ paths, one gets a map $\mathcal D\to X$ where $\mathcal D\defeq R'\cup D\cup (\bigcup_\alpha D_\alpha)$ is such that the composition $\partial_c \mathcal D\to X$ is naturally identified with $\partial N_{A}(v)\subset R$. By \Corref{2x2} applied to $N_{A}(v)$ and $\mathcal D$, there exists a $0$-cell $v'$ in $\mathcal D$ whose image is adjacent to all the neighbors of $v$ in $R$.

    We will show that $v'$ is contained in $D$.  By \Corref{isodiametric}, the diameter of any $D_\alpha$ is at most $3h$ for each $\alpha$. Since $v$ has distance at least $5h$ to $\alpha$, we see that $v$ cannot be adjacent to any $0$-cell in the image of $D_\alpha\to X$. On the other hand, no $0$-cell of $R$ aside from $v$ is adjacent to all the neighbors of $v$ in $R$.  Thus $v'$ must be contained in $D$ and so maps to a $0$-cell of $A'$. The flat $F$ differs by a vertex move from a flat with $v$ replaced by $v'$, since $v'$ has the same neighbors as $v$ in $A$.

    Thus we have shown that given any $0$-cell $v$ of $A$, there exists a vertex move of $F$ replacing it with a $0$-cell of $A'$.  Applying this successively to $0$-cells in an enumeration of the $0$-cells of $A$, we see that $F$ and $F'$ are equivalent.
\end{proof}

\begin{defn}
\defnlabel{thickflat}
    Let $F \colon \flat \to X$ be a flat in a quadric complex. The \emph{thickening of $F$}, denoted $\Th(F)$, is the full subcomplex spanned on all flats at finite Hausdorff distance from $F$. We also call $\Th(F)$ a \emph{thick flat}.
\end{defn}

\begin{thm}
\thmlabel{thickflat}
    Let $F \colon \flat \to X$ be a flat in a quadric complex $X$. There exists a combinatorial map $r:\Th(F)\to \flat$ such that each fibre $r^{-1}(v)$ of a vertex $v\in \flat$ is a set of pairwise nonadjacent vertices, and vertices in fibres $r^{-1}(v)$ and $r^{-1}(w)$ are adjacent if and only if $v$ and $w$ are adjacent. 
\end{thm}

\begin{proof}
    Let $v \in \Th(F)$.  We will show that there is a unique flat $F_v \colon \flat \to X$ whose image contains $v$ such that $F_v$ is obtained from $F$ by a single vertex move.  Since $v \in \Th(F)$ we know that $v$ is in the image of $F'$ for some flat $F'$ equivalent to $F$. By \Lemref{equivalentFlats}, there is a sequence $F=F_0, F_1, F_2, \cdots$ of flats such that $F_{i+1}$ is obtained from $F_i$ by a single vertex move and $F_i \to F'$.  There is a least $i$ for which the image of $F_i$ contains $v$. Suppose $v' = F_{i-1}(u)$ and $v = F_i(u)$. That is, $F_{i-1}$ and $F_i$ differ by a vertex move on $u$. By the last sentence of \Lemref{vertexMovesCombine}, $F$ differs from a flat $F_v$ by a vertex move on $u$, where $v$ is in the image of $F_v$.  Since the flats are isometrically embedded, this $F_v$ is unique.
    
    For $v \in \Th(F)$, we define $r(v)$ to be the unique vertex $v'$ of $\flat$ mapping to $v$ under $F_v$.  Now, consider $v \in r^{-1}(v')$ and $w \in r^{-1}(w')$ for some $v',w' \in \flat$.  If $v' = w'$ then $v$ and $w$ are both adjacent to any neighbour of $v'$ in the image of $F$ and so $v$ and $w$ cannot be adjacent as $X^1$ is bipartite.  If $v' \neq w'$ then, by \Lemref{vertexMovesCombine} there is a flat obtained from $F$ by at most two vertex moves with $v$ and $w$ the images of $v'$ and $w'$, respectively. Thus $v$ and $w$ are adjacent if and only if $v'$ and $w'$ are adjacent.
\end{proof}

\section{The Flat Torus Theorem}
\seclabel{flat_torus_theorem}

\begin{lem}
\lemlabel{normalFlat}
Suppose $G$ acts on a quadric complex $X$, and a normal $\Z^2$ subgroup of $G$ acts metrically properly on a flat $F$ in $X$. Then $G$ stabilizes $\Th(F)$.   
\end{lem}

Note that an action of $\Z^2$ on a flat $F$ is metrically proper if and only if it is proper if and only if it is free if and only if it is cocompact. 

\begin{proof}
    Let $x\in \Th(F)$, and consider a flat $F'$ at finite Hausdorff distance from $F$ containing $x$. For any $g\in G$, since $g\Z^2g^{-1}=\Z^2$ acts cocompactly on $gF$, the flats $F$ and $gF$ are at finite Hausdorff distance. Because $gF$ and $gF'$ are also at finite Hausdorff distance, $gF'$ and $F$ are at finite Hausdorff distance. In particular $gx\in \Th(F)$.
\end{proof}

\begin{thm}
\thmlabel{flatTorus}
    Let $\Z^2$ act metrically properly on a quadric complex $X$.  Then $\Z^2$ acts cocompactly on a $\Z^2$-invariant flat $F$ in $X$.  Moreover, every vertex $x$ in $\Th(F)$ is contained in a cocompact $\Z^2$-invariant flat.
\end{thm}
\begin{proof}
Suppose $\Z^2$ acts metrically properly (hence freely) on a quadric complex $X$ with quotient $\phi:X\to X/\Z^2$. 
By \Apxthmref{surfaceMap} in the appendix, there exists a map from a torus $f:T\to X/\Z^2$ so that $f_*:\pi_1 T\to \pi_1 X/\Z^2$ is an isomorphism. 
Let $\mathsf d$ be the pseudometric on $\widetilde T$ induced by pulling back the metric on $X$ via the $\Z^2$-equivariant map $\tilde f:\widetilde T\to X$.
The action of $\Z^2$ on $(\widetilde T, \mathsf d)$ is metrically proper and is cofinite on $\widetilde T^{(0)}$. 
Thus, letting $x_1,\dots, x_k\in \widetilde T^{(0)}$ be a set of orbit representatives, there are finitely many elements $g\in \Z^2$ that satisfy $\mathsf d(x_i,gx_i)\leq 6$ for some $i\in \{1,\dots, k\}$.
Passing to a finite-index subgroup $H<\Z^2$ that omits all such $g\in \Z^2$, any $h\in H$ acting on $(\widetilde T, \mathsf d)$ has translation length at least $6$.
Indeed, for any $x,y\in \widetilde T$ in the same $\Z^2$-orbit, $\mathsf d(x,gx)=\mathsf d(y,gy), \forall g\in \Z^2$ since $\Z^2$ is abelian.

As in the proof of \cite{Hoda:2017}[Lemma~1.6], we can perform surgeries on $\widetilde T$ in the following ways:
\begin{enumerate}[label=(\Alph*)]
    \item If there exists a vertex $v\in \widetilde T$ with degree less than $4$, then a locally area-reducing surgery can be performed on $N(v)$.
    \item If $T$ is \emph{locally flat}, i.e. the vertices of $T$ all have degree exactly $4$, and $\widetilde T\to X$ contains a knight move, then there exists an area-preserving surgery on the $2\times 1$ subgrid $D$ supporting the knight move that results in $T$ having a vertex of degree $3$ in $\partial D$. 
\end{enumerate}
Because the translation lengths are at least $6$ with respect to $\mathsf d$ (and hence with respect to the graph metric as well), both (A) and (B) surgeries can be made $H$-equivariantly.
Since the (A) and (B) surgeries do not introduce new vertices, translation lengths remain at least $6$ after an $H$-equivariant (A) or (B) surgery.
Furthermore, an $H$-equivariant (B) surgery can always be followed by an $H$-equivariant (A) surgery. 
Because an $H$-equivariant (A) surgery reduces the number of squares in $T=\widetilde T/H$, any sequence of $H$-equivariant surgeries in which each (B) is followed by an (A) must terminate.

By the above, after performing finitely-many $H$-equivariant (A) surgeries, the vertices of $\widetilde T$ all have degree at least $4$. 
Since $\chi(T)=\chi(\widetilde T/\mathbb Z^2)=0$, a combinatorial Gauss-Bonnet argument implies $T$ is locally flat.
If there exists a knight move for $\widetilde T\to X$, then we may perform an $H$-equivariant (B) surgery followed by finitely many $H$-equivariant (A) surgeries, to obtain a locally flat $T$ of lesser area, again by combinatorial Gauss-Bonnet.
Repeating these steps at most finitely many times, we eventually arrive at a locally flat $T$ such that $\widetilde T\to X$ has no knight moves. 

Since $T$ is locally flat, $\widetilde T$ is isomorphic to $\flat$. 
We claim that $\widetilde T\to X$ is locally minimal and thus a flat by \Thmref{localtoflat}. 
Indeed, by \Thmref{no knight and not minimal implies cone}, if $\widetilde T\to X$ were not locally minimal, then there would exist a vertex $w\in X$ adjacent to infinitely many vertices of the image.
However, this contradicts metric properness of the $\mathbb Z^2$ action on $X$. 
Hence $\widetilde T\to X$ is an $H$-equivariant flat. 
Moving forward, we identify $\widetilde T$ with its image in $X$.

% The action of $H$ on $\widetilde T\subset X$ is free and cocompact. Let $g_1H,\ldots, g_nH$ be the cosets of $H<\Z^2$. Note that $g_i\widetilde T$ is stabilized by $g_iHg_i^{-1}=H$, so the collection of flats $\{g_1\widetilde T,\ldots, g_n\widetilde T\}$ is pairwise at finite Hausdorff distance. In particular, $\Th(\widetilde T)=\Th(g_i\widetilde T)$ for any $i=1,\ldots, n$. Furthermore, the collection is $\Z^2$-invariant since $\{g_i\}$ is a set of coset representatives. Hence . 

Since $H$ is a normal subgroup of $\Z^2$ and acts metrically properly on $\widetilde T$, the thick flat $\Th(\widetilde T)$ is $\Z^2$-invariant by \Lemref{normalFlat}. 
Letting $r:\Th(\widetilde T)\to \flat$ be the map from \Thmref{thickflat}, the action of $\Z^2$ on $\Th(\widetilde T)$ is fiber preserving. 
Hence there exists an action of $\Z^2$ on $\flat$ making $r:\Th(\widetilde T)\to \flat$ equivariant. 
By freeness of the $\Z^2$-action, we can equivariantly choose vertices $v\in r^{-1}(v')$ for each $v'\in \flat$, and the chosen vertices span a $\Z^2$-invariant flat $F\subset \Th(\widetilde T)$. 
Since any $v\in \Th(\widetilde T)$ can be chosen, every $v\in \Th(\widetilde T)$ is contained in a cocompact $\Z^2$-invariant flat $F$.
\end{proof}

% Note to future Zach and Nima: we can't equivariantly pick vertices in the $r^{-1}(u)$ so we don't get an action on a flat in the virtually $\Z^2$ case.

\begin{cor}
    Suppose a virtually $\Z^2$ group $G$ acts metrically properly on a quadric complex $X$.  Then $G$ stabilizes $\Th(F)$ for a flat $F$ in $X$.
\end{cor}
\begin{proof}
    Let $H$ be a finite index, normal subgroup of $G$ isomorphic to $\Z^2$. By \Thmref{flatTorus}, $H$ acts properly on a flat $F$ in $X$. By \Lemref{normalFlat}, $G$ acts on $\Th(F)$.
\end{proof}

\begin{defn}
    Let $G$ be a group acting on a metric space $X$. The \emph{min-set} $\Min(g)$ of an element $g\in G$ is the set of $x\in X$ such that $d(x,gx)\leq d(y,gy)$ for all $y\in X$. The \emph{min-set} of $G$ is the intersection $\Min(G)=\cap_{g\in G}\Min(g)$.
\end{defn}

\begin{thm}
    Let $\Z^2$ act freely on a quadric complex $X$. Let $F$ be a $\Z^2$-invariant flat in $X$. Then $\Min(\Z^2)=\Th(F)$.
\end{thm}

Note for any free action of $\Z^2$ on $X$ there exists a $\Z^2$-invariant cocompact flat $F$ in $X$ by \Thmref{flatTorus}.

\begin{proof}
$\Th(F)\subset \Min \Z^2$: Suppose there exists $x\in \Th(F)$ such that $x\not\in \Min\Z^2$. Then there exist $g\in \Z^2$, $m\in X$ so that $d(m,gm)<d(x,gx)$. For any $n>0$, we have $d(x,g^nx)\leq d(x,m) + d(m,g^nm) + d(g^nm,g^nx)\leq 2d(m,x) + nd(m,gm)$. By \Thmref{flatTorus} $x$ belongs to a $\Z^2$-invariant flat $F$, so $nd(x,gx)=d(x,g^nx)\leq 2d(m,x)+nd(m,gm)$. For $n\gg 0$, we get a contradiction.  

$\Min \Z^2\subset \Th(F)$: Let $m\in \Min\Z^2$ be a minimally translated element, and fix an arbitrary $x\in F$. 
Since $x\in \Min\Z^2$, we have that $d(gx,hx)=d(gm,hm)$ for all $g,h\in \Z^2$. Let $\ell_1$ and $\ell_2$ be the vertical and horizontal lines in $F$ passing through $x$. By cocompactness of the $\Z^2$ action on $F$, there exist nontrivial $g,h\in \Z^2$ so that $g(x)\in \ell_1$ and $h(x)\in \ell_2$. Take geodesics $\gamma$ and $\eta$ joining $m$ to $gm$ and $m$ to $hm$, and let $D$ be a minimal area diagram filling $\partial D\defeq \gamma(g\eta)(h\gamma)^{-1}\eta^{-1}$. We call these four subpaths of $\partial D$ the \emph{sides} of $D$. Their four junctions are the \emph{corners} of $D$. Since $g$, $h$, and $gh$ translate $x$ and $m$ the same distances, the concatenation of two sides of $D$ meeting at a corner is a geodesic. 
Thus opposing corners of $D$ are at distance $|\gamma|+|\eta|$.
Consequently, opposite sides of $D$ are disjoint. 
Indeed, if opposite sides of $D$ intersect, then there would exist a path in the union of the intersecting sides of length less than $|\gamma|+|\eta|$ joining opposite corners of $D$, a contradiction.
Consequently, $D$ is non-singular, i.e. homeomorphic to a disk.

We will prove that $D$ is isomorphic to the product of its sides $\gamma$ and $\eta$. It suffices to show the curvature at the corners $m$, $gm$, $hm$, and $ghm$ in $\partial D$ is $\frac14$, and the curvature at any other vertex of $D$ is zero. 
An inductive argument using the curvature conditions and non-singularity of $D$ imply $D$ is the product of $\gamma$ and $\eta$. 

% \note{maybe cite something}

A negatively curved vertex has curvature at most $-\frac14$, and since the concatenation of adjacent sides are geodesic, no vertex of $\partial D$ can have curvature greater than $\frac14$. 
Since $D$ is non-singular, any two vertices of positive curvature in the interior of a geodesic in $\partial D$ must be separated by a vertex of negative curvature. 
Consequently, the sum of curvatures of the interior vertices of a geodesic in $\partial D$ cannot exceed $\frac14$.
By \Lemref{minareaCAT0}, the interior vertices of $D$ have nonpositive curvature, since $D$ is minimal area. 
By \Lemref{gauss_bonnet}, the sum of curvatures of vertices in $D$ must total one. 
Consider the decomposition $\partial D$ into the two geodesics $\gamma(g\eta)$ and $(h\gamma)^{-1}\eta^{-1}$. 
The interiors of each geodesic can contribute at most $\frac14$ to the total curvature of $D$, thus vertices $m$ and $ghm$ must together contribute $\frac12$ to the total curvature. 
Hence each contributes $\frac14$, by the observation above. 
A symmetric argument shows $gm$ and $hm$ have curvature $\frac{1}{4}$. 

It remains to prove any non-corner vertex of $D$ has curvature zero. We first consider boundary vertices of $D$. The vertex $gm$ in the geodesic $\gamma(g\eta)$ has curvature $\frac14$. As observed above, any pair of positively curved vertices in the interior of $\gamma(g\eta)$ must be separated by at least one vertex of negative curvature. In particular, if there exists a vertex of non-zero curvature in $\gamma$, then the closest such vertex to $gm$ is negatively curved. Similarly, for any corner, a closest vertex of non-zero curvature in the interior of an incident side is negatively curved. Thus, since pairs of positive vertices are separated by at least one negative vertex, either the interior vertices of a side all have zero curvature, or the sum of curvatures is negative. Consequently, if any vertex of $\partial D$ apart from the corners $m$, $gm$, $hm$, and $ghm$ has nonzero curvature, then the total curvature of $\partial D$ is less than one. The interior of $D$ has everywhere nonpositive curvature. Hence, since $D$ has total curvature $1$, every vertex apart from the corners has zero curvature. 

We have concluded that the diagram $D\to X$ is isomorphic to the product of its sides $\gamma$ and $\eta$. Then, by gluing together $\langle g, h\rangle$-translates of $D \to X$, we obtain a $\langle g, h\rangle$-equivariant map $\flat\to X$.  For any pair of vertices $u,v \in \flat$, there exists a geodesic $\gamma$ of $\flat$ containing $u$ and $v$ whose endpoints lie in the orbit $\langle g,h \rangle m$.  Since the restriction of $\flat \to X$ to $\langle g,h \rangle m$ is an isometric embedding, it follows that $\flat \to X$ is isometric and hence is a flat.  As each $\langle g,h \rangle$-translate of $m$ lies within a uniform distance to the corresponding translate of $x$ in $F$, the image of $\flat \to X$ is within finite Hausdorff distance to $F$ and so is contained in $\Th(F)$. 
\end{proof}

\begin{defn}
    An action of a group $G$ on a metric space $X$ is \emph{metrically proper} if for every closed ball $B\subset X$, the set of elements $\{g\in G:gB\cap B\neq \emptyset\}$ is finite.   
\end{defn}

\begin{thm}
    Let $\Z^n$ act metrically properly on a quadric complex $X$. Then $n\leq 2$. In particular, if $\Z^n$ acts freely on a locally finite quadric complex $X$, then $n\leq 2$. 
\end{thm}

\begin{proof}
    Suppose $n>2$, and consider a $\Z^2$ subgroup of $\Z^n$. By \Thmref{flatTorus}, $\Z^2$ acts on a flat $F$ in $X$. By \Lemref{normalFlat}, $\Z^n$ acts on $\Th(F)$. The map $r:\Th(F)\to \flat$ from \Thmref{thickflat} induces an action of $\Z^n$ on $\flat$. 
    
    Fix a $g\in \Z^n\setminus \Z^2$. Up to passing to a power, $g$ acts by translation on $\flat$ and there exists $h\in \Z^2$ so that $gh$ acts trivially on $\flat$. Thus $gh\neq 1$ acts on a fiber of $r$, contradicting metric properness of the action.
\end{proof}

\appendix 

\section*{Appendix. Combinatorial surface maps}
\seclabel{surface_maps}

\begin{apxdefn}
	Let $X$ be a 2-complex.
    The \emph{disconnection degree} $\ddeg x$ of a 0-cell $x\in X$ is the number of components in $\lk x$. 
    The 0-cell $x$ is \emph{singular} if $\lk x$ is a single vertex or disconnected. 
    A 1-cell of $X$ is \emph{singular} if it does not belong to the boundary of any 2-cell. 
\end{apxdefn}

\begin{apxdefn}\apxdefnlabel{}
    A \emph{degenerate map} of combinatorial $2$-complexes is a continuous map $f:Y\to X$ such that for each open cell $C$ of $Y$
    \begin{itemize}
        \item $f(C)$ is an open cell of $X$ and $\dim f(C)\leq \dim C$,
        \item if $\dim C=\dim f(C)$ for an open cell $C$ of $Y$, then $f|_C$ is a homeomorphism onto $f(C)$,
        \item if $\dim C=2$ and $\dim f(C)=1$, then 
        there are exactly two open $1$-cells of $\partial_c C$ that map homeomorphically onto the open $1$-cell $f(C)$ and any other $1$-cell of $\partial_cC$ maps to a $0$-cell of $\partial f(C)$, and
        \item if $\dim C=1$ and $\dim f(C)=0$, then $C$ is a $1$-cell of some $2$-cell $C'$ such that $\dim f(C')=0$.
    \end{itemize}
    A \emph{thick $0$-cell} (resp. \emph{a thick $1$-cell}) is a $2$-cell $C$ of $Y$ such that $\dim f(C)=0$ (resp. $\dim f(C)=1$).
\end{apxdefn}

\begin{apxlem}\apxlemlabel{thickened_diagrams}
    Let $X$ be a combinatorial $2$-complex, and let $S$ be a closed, orientable surface of nonzero genus. 
    If $f_*:\pi_1S\to \pi_1X$ is an embedding, then there exists a combinatorial structure on $S$ and a degenerate map $f:S\to X$ such that $f$ induces~$f_*$. 
\end{apxlem}

\begin{figure}[h]
        \centering
        \includegraphics[width=0.5\textwidth]{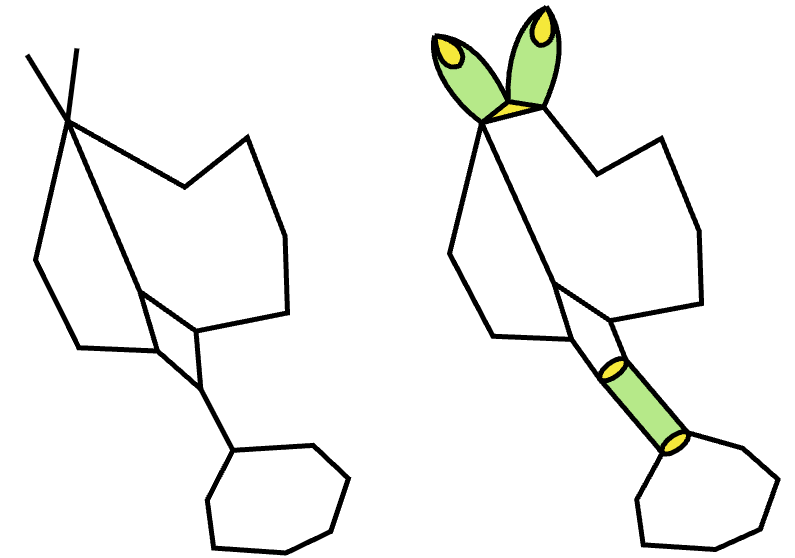}
        \caption{An example of a singular diagram $D$ (left) and its thickening $\mathbf D$ (right). The blown-up singular edges are colored green, and the blown-up singular vertices are colored yellow. Note the degree one vertices in $D$ become $1$-gons in $\mathbf D$.}
        \label{fig:thickDiagram}
\end{figure}

\begin{proof}
    Consider the presentation $\pi_1 S=\langle a_1,b_1,\ldots, a_g,b_g \mid [a_1,b_1]\cdots[a_g,b_g] \rangle$. 
   The generators $\{a_i,b_i\}_i$ in the presentation correspond to closed paths $\{\alpha_i,\beta_i\}_i$ in $X$. 
   Concatenating these closed paths as in the word $[\alpha_1,\beta_1]\ldots[\alpha_g,\beta_g]$ defines a nullhomotopic, closed path $P\to X$. 
   Let $D\to X$ be a diagram in $X$ with $\partial_c D\to D\to X$ equivalent to $P\to X$.
   The decomposition $[\alpha_1,\beta_1]\cdots[\alpha_g,\beta_g]$ of $\partial_c D$ defines a quotient $\overline D$ of $D$ and map $\overline D \to X$ as follows: for each $\gamma\in \{\alpha_i,\beta_i\}_i$, we identify the two subpaths of $D$ equivalent to $\gamma$ in $X$. 
    
    We define a ``thickening'' $\mathbf D$ of $D$ as follows: Let $\{x_i\}$ and $\{e_j\}$ be the singular 0-cells and 1-cells of $D$. 
    Let $\{D_k\}$ be the set of closures of the components of $D-\big(\bigcup_i x_i \cup \bigcup_j \Int(e_j)\big)$. Note that each $D_k$ is homeomorphic to a closed $2$-ball. 
    Let $y$ be 0-cell of $D_k$ that is singular in $D$. Pick an arbitrary corner of a 2-cell $R$ of $D_k$ containing $y$ and blow up that corner to contain a $1$-cell $\mathbf y$ in place of $y$. 
    Repeating this procedure for each singular 0-cell $y$ in $D$ and each $D_k$ containing $y$, we obtain complexes $\{\mathbf D_k\}_k$. 
    Each $\mathbf D_k$ is still homeomorphic to a closed $2$-ball. For each $k$, there is a quotient map $\mathbf{D}_k \to D_k$ collapsing each $\mathbf{y} \subseteq \mathbf{D}_k$ onto $y$.  
    For each $j$, let $\mathbf e_j$ be the rectangle $e_j\times I$. 
    For each $i$, let $\mathbf C_i$ be a $\ddeg {x_i}$-gon. Note we allow $\partial\mathbf C_i$ to be a loop (i.e. a 1-gon), which happens when $x_i$ has both disconnection degree and degree equal one. 
    We identify pairs of the new blown-up 1-cells in $\mathbf D_k$, $\mathbf e_j$ to the $\mathbf C_i$ respecting the cylic order of the $D_k$, $e_j$ that intersect $x_i$ in $D$. 
    See Figure~\ref{fig:thickDiagram}. 
    Each 1-cell of a $\mathbf C_i$ is identified with either a blown-up 1-cell $\mathbf y$ in some $\mathbf D_k$ or one of the two 1-cells in $\partial e_i\times I$ in some $\mathbf e_i$. 
    We let $\mathbf D$ denote the complex obtained by making all such identifications. 
    Note $\mathbf D$ is homeomorphic to a closed $2$-ball. 

    Via the equivalence between $\partial_c D$ and $\partial_c\mathbf D$ there is a decomposition of $\partial_c\mathbf D$ as $[\alpha_1,\beta_1]\cdots[\alpha_g,\beta_g]$, and identifying pairs of equivalent paths, we obtain a quotient $\mathbf D \to \overline {\mathbf D}$.
    There exists a degenerate map $\mathbf D\to D$ sending each $\mathbf C_i$ to $x_i$ and $\mathbf e_j$ to $e_j$ (by collapsing the $I$ direction of $\mathbf e_j$) that descends to a degenerate map $\overline{\mathbf D}\to \overline D$.
    The composition $f:\overline{\mathbf D}\to\overline D\to X$  induces $f_*$. 
    Note the thick 0-cells and thick 1-cells of $\overline{\mathbf D}$ are exactly $\{\mathbf C_i\}_i$ and $\{\mathbf e_j\}_j$.
    Since $\mathbf D$ is homeomorphic to a closed $2$-ball, the quotient $\overline{\mathbf D}$ is homeomorphic to $S$. 
\end{proof}

\begin{apxthm}
\apxthmlabel{surfaceMap}
    Let $X$ be a combinatorial 2-complex and $S$ a closed, orientable surface of nonzero genus. Let $\phi_*:\pi_1 S\to \pi_1 X$ be an embedding. Then there exists a combinatorial structure on $S$ and a combinatorial map $\phi:S\to X$ inducing $\phi_*$. 	
\end{apxthm}

\begin{proof}

% onto $\phi_*(\pi_1S) < \pi_1 X$. In fact, since surface groups are Hopfian \cite{}, $q_*$ maps

	By \Apxlemref{thickened_diagrams}, there exists a combinatorial structure on $S$ and a degenerate map $\phi:S\to X$ inducing $\phi_*$.
    Choosing $\phi$, $S$ so that $S$ has a minimal number of cells, we claim $\phi$ is a combinatorial map.  
    Let $\{\mathbf C_i\}_i$ and $\{\mathbf e_j\}_j$ be the set of thick $0$-cells and thick $1$-cells, respectively.

    We will prove that the set of thick $0$-cells is empty. 
    If there exists a thick $0$-cell $\mathbf C_i$, then there exists a simple closed path $P$ in $\mathbf C_i \subset S$. 
    Since $\phi(P)$ is a point and $\phi_*$ is injective, $P$ must bound an embedded disk $D$ in $S$. 
    The quotient map $S\to S/D$ is a homotopy equivalence, the quotient $S/D$ is homeomorphic to $S$, and $\phi$ descends to a degenerate map on $S/D$, contradicting minimality.

    Since there are no thick $0$-cells, each thick $1$-cell $\mathbf e_j$ is a bigon, and $\phi(\mathbf e_j)$ is an edge of $X$, where $\phi\vert_{\mathbf e_j}$ factors through a retraction of $\mathbf e_j$ onto one of its edges.
    We claim there are no thick $1$-cells. 
    Indeed, if $\mathbf e_j$ is a thick $1$-cell, then its edges are distinct, since otherwise $\mathbf e_j\subset S$ would be homeomorphic to a sphere, a contradiction.
    Thus deleting $\Int(\mathbf e_j)$ and identifying the edges of $\partial\mathbf e_j$, we obtain a surface with fewer cells, contradicting minimality.
    
    We conclude by noting that a degenerate map with no thick $0$-cells or $1$-cells is a combinatorial map.
\end{proof}

\begin{apxdefn}[Reduced, cancellable pair]
    Let $f \colon Y\to X$ be a combinatorial map of combinatorial $2$-complexes and let $R$ and $R'$ be $2$-cells of $Y$.
    The pair $R,R'$ is \emph{cancelable} along an oriented edge $e$ of $Y$ if the attaching maps of $R$ and $R'$ decompose as distinct paths $e e_1 e_2 \cdots e_k$ and $e e'_1 e'_2 \cdots e'_k$ in $Y$ that become identical paths in $X$ when composed with $f$.
    The combinatorial map $Y\to X$ is \emph{reduced} if it has no cancellable pairs.
\end{apxdefn}

If $Y\to X$ is reduced, then it is locally injective away from the $0$-cells of $Y$.

% \begin{apxdefn}[Reduced, cancellable pair]
%     Let $Y\to X$ be a combinatorial map of combinatorial $2$-complexes. 
%     Two $2$-cells $R,R'$ of $Y$ with attaching maps $\partial_c R \to Y$ and $\partial_c R' \to Y$ form a \emph{cancellable pair} if there exist 
%     \begin{itemize}
%         \item combinatorial maps $e\to \partial_c R$ and $e\to \partial_c R'$ from a $1$-cell $e$ and
%         \item an isomorphism $\partial_cR'\to \partial_c R$ 
%     \end{itemize} 
%     making the following diagram commute

%     \[ \begin{tikzcd}[column sep=tiny, row sep = tiny]
%         e \ar[rr] \ar[ddd] & \  & \partial_c R \ar[ddd] \ar[ddddr, bend left=15] \\
%         \ \\
%         \ \\
%         \partial_c R' \ar[rr] \ar[uuurr] \ar[drrr, bend right=15] & & Y \ar[dr] \\
%         & & & X.
%     \end{tikzcd} \]
%     The map $Y\to X$ is \emph{reduced} if there exist no cancellable pairs.
% \end{apxdefn}

\begin{apxdefn}[Surface diagrams]
    A \emph{surface diagram} is an embedding of a finite $2$-complex in a closed surface $D\to S$ such that the embedding induces a cell structure on $S$ where $S-D$ is either empty or a single $2$-cell denoted $C_\infty$. 
    The \emph{boundary cycle} $\partial_cD\to D$ is the attaching map of $C_\infty$.
\end{apxdefn}

Note a disk diagram is a surface diagram where $S$ is the $2$-sphere and $D-S$ is a $2$-cell.

\begin{apxthm}
\apxthmlabel{surfaceReduction}
	Let $X$ be a combinatorial 2-complex and $\phi:S\to X$ a $\pi_1$-injective combinatorial map from a closed, orientable, combinatorial surface of nonzero genus to $X$. 
    If $\phi$ is not reduced, then there exists a combinatorial structure $\overline S$ on $S$ with fewer cells and a $\pi_1$-injective combinatorial map $\bar\phi:\overline S\to X$ with $\phi_*(\pi_1S)=\bar\phi_*(\pi_1\overline S)$. 
\end{apxthm}

\begin{proof}
	Suppose $R$, $R'$ are a cancelable pair in $S$ meeting at an edge $e$. 
    Then the attaching maps of $R$ and $R'$ can be expressed as closed paths beginning with $e$ that are distinct in $S$ but become identical in $X$ when composed with $\phi$.
    If $R$ and $R'$ are distinct $2$-cells of $S$, remove the open cells $R$, $R'$, and $e$ from $S$ to obtain a complex $S_1$. 
    If $R$ and $R'$ are the same $2$-cell, then because $S$ is orientable, the union of open cells $R\cup e$ is an open anulus and in this case we remove only the open cell $R$ from $S$ to obtain a complex $S_1$.
    In both cases, the inclusion $S_1\hookrightarrow S$ is a surface diagram.
    Moreover, $S - S_1$ a $2$-cell $C_\infty$ whose attaching map can be expressed as a path $\partial_cC_\infty = e^{-1}_k\dots e^{-1}_1e'_1\dots e'_k$ such that $e_1 \dots e_k$ and $e'_1\dots e'_k$ are distinct paths in $S_1$ that become identical in $X$ when composed with $\phi$.
    Let $\phi_1:S_1\to X$ be the restriction of $\phi$ to $S_1$, and note that $(\phi_1)_*$ surjects onto $\phi_*(\pi_1S)$. We iteratively construct a sequence of surface diagrams $S_i\subset S$ and continuous (potentially non-combinatorial) maps $\phi_i:S\to X$, $i=1,2,\dots$, such that 
    \begin{itemize}
        \item $\phi_i|_{S_i}$ is combinatorial,
        \item $\phi_i$ induces a $\pi_1$-surjection onto $\phi_*(\pi_1S)$, 
        \item each $S_{i+1}$ is either a subspace, quotient, or quotient of a subspace of $S_i$, and $\phi_{i+1}|_{S_{i+1}}$ is the natural restriction, descent, or descent of a restriction of $\phi_i$, 
        % \item $e_i, e'_i$ remain as edges of $S_i$ and have the same initial vertex, 
        % \item $e_i,\dots,e_k, e'_i$ have either been identified in $S_{i+1}$ or no longer remain in $S_{i+1}$, and
        \item if $S-S_i$ is a $2$-cell $C_\infty$, then $\partial_cC_\infty=e^{-1}_k\dots e^{-1}_ie'_i\dots e'_k$, where the previous point allows us to make sense of $e_j,e'_j$ as potential edges of $S_i$.
    \end{itemize}
    Note the $e_i, e'_i$ are oriented edges, so we may refer to their initial and terminal vertices.
    Also, it follows from the above assumptions that $(\phi_i)_*$ is $\pi_1$-injective since $\pi_1S$ is Hopfian.

    At each step, $S_{i+1}$ is obtained from $S_i$ either by identifying $e_i$ and $e'_i$, by removing $e_i = e'_i$, or by removing an open disk from $S_i$ and identifying $e_i$ and $e'_i$.  
    We now describe this process in detail below.

    \begin{figure}[h]
    \centering
    \begin{subfigure}[h]{0.65\textwidth}
        \centering
        \includegraphics[width=\textwidth]{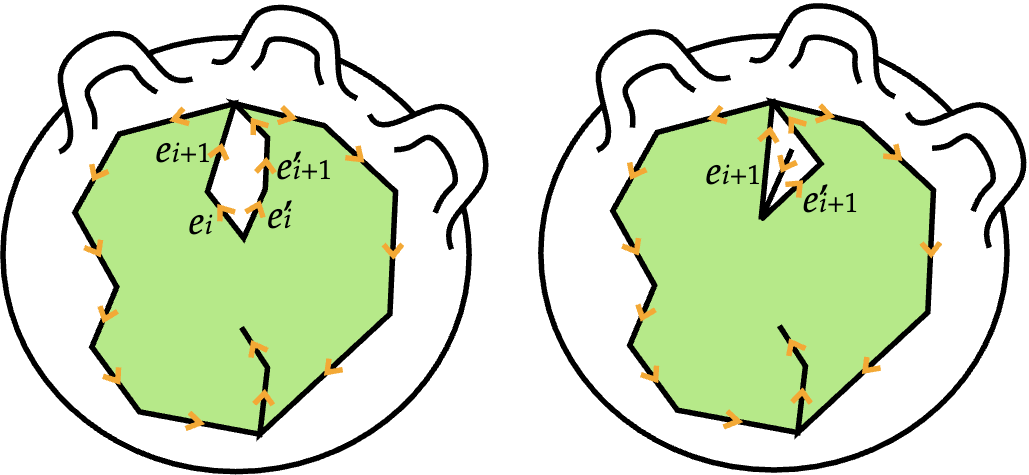}
        \caption{Collapsing a corner of $C_\infty$ when $e_i, e'_i$ share a single endpoint.}	
        \label{fig:surgery1}
    \end{subfigure}
	\ \ \ \ \ 
    \begin{subfigure}[h]{0.65\textwidth}
        \centering
        \includegraphics[width=\textwidth]{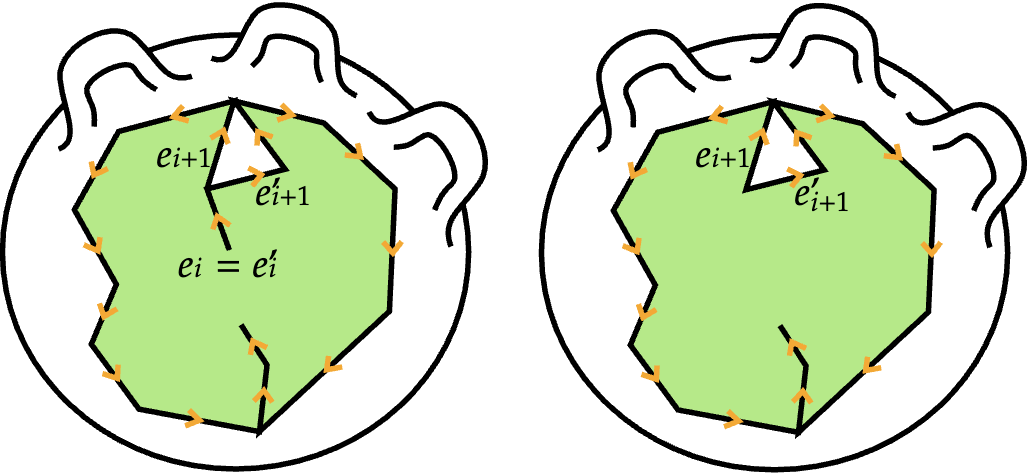}
        \caption{Removing $e_i=e'_i$.}	
        \label{fig:surgery2}
    \end{subfigure}
	\ \ \ \ \ 
    \begin{subfigure}[h]{0.65\textwidth}
        \centering	
	\includegraphics[width=\textwidth]{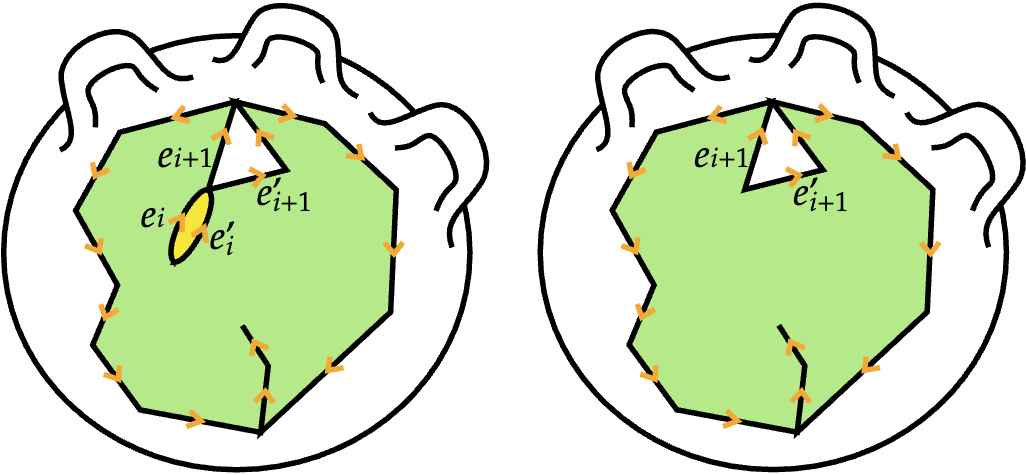}
        \caption{$e^{-1}_ie'_i$ bounds a disk $D$, colored yellow.}	
        \label{fig:surgery3}
    \end{subfigure}
    \caption{Various ways of passing from $S_i$ to $S_{i+1}$. The $2$-cell $C_\infty$ is colored green, and the oriented edges $e_i,\dots, e_k$ and $e'_i,\dots e'_k$ are indicated with orange arrows.}
    \label{fig:surgery}
\end{figure}
    
    Suppose we have constructed $\phi_i:S\to X$. 
    If $S_i=S$, then we can set $\bar S = S_i$ and $\bar \phi = \phi_i$ and we are done, so we may assume that $S-S_i$ is a $2$-cell $C_\infty$.
    If $e_i, e'_i$ share a single endpoint, then $e_i,e'_i$ form a corner of $C_\infty$. 
    We collapse this corner, identifying $e_i$ and $e'_i$, to form $S_{i+1}$ in a quotient of $S$ that we can identify with $S$ via a homeomorphism.
    See Figure~\ref{fig:surgery1}.
    We then let $\phi_{i+1}$ be the descent of $\phi_i$ to the quotient. 

    If $e_i,e'_i$ share both endpoints, then $e^{-1}_ie'_i$ is a nullhomotopic path in $S$ since $\phi_i(e_i)=\phi_i(e'_i)$ and $(\phi_i)_*$ is $\pi_1$-injective.
    Hence either $e_i=e'_i$ or $e^{-1}_ie'_i$ bounds a disk in $S$. 
    If $e_i=e'_i$, then the shared intial vertex $v$ has degree one, so we let $S_{i+1}=S_i-(\interior{e_i}\cup v)$ and let $\phi_{i+1}=\phi_i$.
    See Figure~\ref{fig:surgery2}.
    If $e^{-1}_ie'_i$ bounds a disk $D\subset S$ then there are two cases: either $C_{\infty} \subset D$ or $C_{\infty} \subset S \setminus \interior D$.  
    In either case, removing all open cells of $S_i$ contained in $\interior D$ and collapsing the bigon $e^{-1}_ie'_i$ along with $D$ onto a single edge (identifying $e_i$ and $e'_i$), we obtain $S_{i+1}$ in a quotient of $S$ identifiable with $S$ as above.
    We let $\phi_{i+1}$ be the descent of $\phi_i|_{S\setminus \interior D}$ to the quotient.
    Figure~\ref{fig:surgery3} illustrates the case where $C_{\infty} \subset S \setminus \interior D$.
    Note that in the case $C_\infty\subset D$ we have $S_{i+1}=S$, so as described above, the sequence terminates in the following step.

    Going from $S_i$ to $S_{i+1}$, the length of $\partial_cC_\infty$ always decreases, so this process must eventually terminate.
 \end{proof}

The above theorems have the following consequence: 
If $X$ is a combinatorial $2$-complex with a surface subgroup $H<\pi_1X$, then by \Apxthmref{surfaceMap} there exists a combinatorial map $\phi:S\to X$ so that $\phi_*:\pi_1S\to H$ is an isomorphism. 
If $\phi:S\to X$ is chosen so that $S$ has the minimal number of $2$-cells, then $\phi$ is reduced by \Apxthmref{surfaceReduction}. 
Thus inclusions of surface subgroups of combinatorial $2$-complexes are induced by reduced maps from combinatorial surfaces.

\bibliographystyle{plain}
\bibliography{nima,more}
\end{document}